# OPTIMAL REINSURANCE/INVESTMENT PROBLEMS FOR GENERAL INSURANCE MODELS

By Yuping Liu and Jin Ma[1]

*Wachovia Corporation and University of Southern California*

In this paper the utility optimization problem for a general insurance model is studied. The reserve process of the insurance company is described by a stochastic differential equation driven by a Brownian motion and a Poisson random measure, representing the randomness from the financial market and the insurance claims, respectively. The random safety loading and stochastic interest rates are allowed in the model so that the reserve process is non-Markovian in general. The insurance company can manage the reserves through both portfolios of the investment and a reinsurance policy to optimize a certain utility function, defined in a generic way. The main feature of the problem lies in the intrinsic constraint on the part of reinsurance policy, which is only proportional to the claim-size instead of the current level of reserve, and hence it is quite different from the optimal investment/consumption problem with constraints in finance. Necessary and sufficient conditions for both well posedness and solvability will be given by modifying the "duality method" in finance and with the help of the solvability of a special type of *backward stochastic differential equations*.

**1. Introduction.** Optimization proportional reinsurance problems have been considered by many authors in recent years. We refer to the books of Gerber [6] and Bühlmann [3] for the basic idea of proportional reinsurance, and to, say, [7] for the treatment for diffusion models. However, in most of the previous works the dynamics of the reinsurance problems, that is, the reserve processes, were usually restricted to the rather simplistic model, such as classical Cramér–Lundberg model or its simple perturbations. One of the consequences of such settings was that the results and methodology used

Received November 2006; revised November 2008.
[1]Supported in part by NSF Grants 05-05427 and 08-06017.

in solving such problems depend, explicitly or implicitly, on the Markovian nature of the reserve process. The generalization of these results to a more realistic environment is therefore rather difficult. In fact, it is still rather fresh.

This paper is an attempt in this direction. We shall consider a generalized insurance model as was proposed in Ma–Sun [14]. More precisely, let us consider a risk reserve process, denoted by $X_t$, that takes the following form:

$$(1.1) \qquad X_t = x + \int_0^t c_s(1+\rho_s)\,ds - \int_0^{t+}\int_{\mathbb{R}_+} f(s,z) N_p(ds\,dz),$$

where $c_s$ is the premium rate process and $\rho = \{\rho_t\}_{t\geq 0}$ is the so-called safety loading process. The last stochastic integral represents a general claim process in which $N_p$ is the counting measure generated by stationary Poisson point process $p$; and $f$ represents the intensity of the jumps (detailed characterizations of these quantities will be given in Section 2).

Our optimization problem is based on the following consideration: we suppose that the insurance company can manage its reserve, whence risk, in three ways: investment, (proportional) reinsurance and consumption. More precisely, we assume that the insurance company puts its reserve in a financial market that contains 1 riskless account and some risky assets, and it is allowed to change its investment positions continuously. Also, we assume that the insurance company can divert (cede) a fraction of the incoming claims, while yielding a fraction of its premium at the same time, to a reinsurance company. Finally, the insurance company is also allowed to "consume" (in the form of "dividend," "refund," etc.). The goal of the insurance company is then to optimize certain utility by managing the investment portfolio, reinsurance policy and the consumption. We should note that since a reinsurance policy must take values in $[0,1]$, our optimization problem seems to resemble the utility optimization problem with portfolio constraints. We refer the readers to, for example, Karatzas–Shreve [10, 11] for the optimal investment/consumption problems with continuous models, to Xue [20] for their jump-diffusion counterparts and to Cvitanic–Karatzas [4] for the results involving portfolio constraints.

In order to avoid over-complicating the model, in this paper we allow a relaxed "admissibility" on the portfolio, namely, we shall allow short selling and borrowing with the same interest rate so that no restrictions are needed on the bounds and signs of the portfolios. The main feature of our problem, compared to the utility optimization problems in finance, comes from the nature of insurance and reinsurance. For example, in an insurance model the jumps come from the claims, which is independent of the financial market, hence not hedgeable. Furthermore, since it is not practical to assume that the reinsurance policy is proportional to the current level of reserve, the dynamics of our "wealth" (reserve) process cannot be formulated as a homogeneous linear SDE, hence keeping the reserve from "ruin" throughout a given duration is by no means trivial. As a matter of fact, such a nonhomogeneity, along with the constraint on the size of reinsurance control (as a proportion), causes the main technical difficulties in this work (compared to, e.g., [20]). To our best knowledge, such problems have not been fully explored, especially under an actuarial context, and it has sufficient novelty even in the utility optimization literature.

The main results of this paper focus on two aspects: the well posedness of the optimization problem and the actual resolution of the optimal strategy. The first part of the results include the study of admissible strategies, and the actual existence of such strategies. After a careful study of the reinsurance structure, via the so-called "*profit-margin principle*," we derive a reasonable risk reserve model with reinsurance and investment. Such a model is a natural extension of the simplest ones one usually sees in the elementary actuarial literature (without diffusion approximations). The admissibility of the portfolio/reinsurance/consumption triplet is then defined so that the insurance company does not go default over a given planning horizon. Due to the constraint on the reinsurance part, the existence of such admissible triplet becomes a rather technical issue. In fact, the verification of the existence of admissible strategy relies on a new result on the so-called *backward stochastic differential equations*, which is interesting in its own right. Our main result on the existence of admissible strategies is then proved along the lines of the result of [5], with some necessary modifications. We should note that our reinsurance policy has to depend on the sizes of the claims. Technically, such a restriction can be removed if the process $S_t$ has fixed size jumps (cf., e.g., [20]), but this is not of significant interest because it will exclude even the simplest compound Poisson claim processes.

Finally, we would like to point out that our utility optimization problem are formulated slightly different from the traditional ones, due to some technical assumptions needed in order to guarantee the existence of admissible strategies. In particular, we will require that the utility function for the terminal reserve to be a "truncated" version so that the terminal wealth of the optimal reserve is bounded. We should note that every utility function can be approximated by a truncated sequences, thus an "$\varepsilon$-optimal" strategy could be produced using our result. Also, it is worth noting that our final result rely on the solvability of a special "forward–backward stochastic differential equation" (FBSDEs for short; cf., e.g., Ma–Yong [15] for more details on such equations). But the existence of the solution to the present FBSDE is by no means trivial, and seems to be beyond the scope of all the existing results. We will not pursue all these issues in this paper due to the length of the paper, but we hope to be able to address them in our future publications.

This paper is organized as follows. In Section 2 we give the necessary preliminaries about our model. In Section 3 we describe the admissibility of the investment–reinsurance–consumption strategies and introduce some equivalent probability measures that are important in our discussion. In Section 4 we introduce the wider-sense admissible strategies and prove the existence of such strategies, and in Section 5 we derive a sufficient condition for the existence of true admissible strategies. The last section is devoted to the utility optimization problem.

**2. Preliminaries and reserve model formulations.** Throughout this paper we assume that all uncertainties come from a common complete probability space $(\Omega, \mathcal{F}, \mathbb{P})$ on which is defined a $d$-dimensional Brownian motion $W = \{W_t : t \geq 0\}$ and a stationary Poisson point process $p$. We assume that $W$ and $p$ are independent, which will represent the randomness from the financial market and the insurance claims, respectively. For notational clarity, we denote $\mathbf{F}^W = \{\mathcal{F}_t^W : t \geq 0\}$ and $\mathbf{F}^p \stackrel{\triangle}{=} \{\mathcal{F}_t^p : t \geq 0\}$ to be the filtrations generated by $W$ and $p$, respectively, and denote $\mathbf{F} = \mathbf{F}^W \otimes \mathbf{F}^p$, with the usual $\mathbb{P}$-augmentation such that it satisfies the *usual hypotheses* (cf., e.g., Protter [16]). Furthermore, we shall assume that the point process $p$ is of class (QL) (cf. [8] or [9]), and denote its corresponding counting measure by $N_p(dt\,dz)$. The compensator of $N_p(dt\,dz)$ is then $\hat{N}_p(dt\,dz) = E(N_p(dt\,dz)) = \nu(dz)\,dt$, where $\nu(dz)$ is the Lévy measure of $p$, satisfying $\nu(\mathbb{R}_+) < \infty$, where $\mathbb{R}_+ \stackrel{\triangle}{=} (0, \infty)$.

Let us specify some notation in this paper. Let $\mathbb{E}$ be a generic Euclidean space. Regardless of its dimension we denote $\langle \cdot, \cdot \rangle$ and $|\cdot|$ to be its inner product and norm, respectively. The following spaces will be frequently used:

- $C([0,T]; \mathbb{E})$ is the space of all $\mathbb{E}$-valued continuous functions on $[0,T]$;
- for any sub-$\sigma$-field $\mathcal{G} \subseteq \mathcal{F}_T$ and $1 \leq p < \infty$, $L^p(\mathcal{G}; \mathbb{E})$ denotes the space of all $\mathbb{E}$-valued, $\mathcal{G}$-measurable random variables $\xi$ such that $E|\xi|^p < \infty$. As usual, $\xi \in L^\infty(\mathcal{G}; \mathbb{E})$ means that it is $\mathcal{G}$-measurable and bounded;
- for $1 \leq p < \infty$, $L^p(\mathbf{F}, [0,T]; \mathbb{E})$ denotes the space of all $\mathbb{E}$-valued, $\mathbf{F}$-progressively measurable processes $\xi$ satisfying $E \int_0^T |\xi_t|^p\,dt < \infty$. The meaning of $L^\infty(\mathbf{F}, [0,T]; \mathbb{E})$ is defined similarly;
- $F_p$ (resp., $F_p^2$) denotes the class of all random fields $\varphi : \mathbb{R}_+ \times \mathbb{R}_+ \times \Omega \to \mathbb{R}_+$, such that for fixed $z$, the mapping $(t, \omega) \mapsto \varphi(t, z, \omega)$ is $\mathbf{F}^p$-predictable, and that

$$E \int_0^T \int_{\mathbb{R}_+} |\varphi(s,z)| \nu(dz)\,ds < \infty$$

(2.1)
$$\left(\text{resp. } E \int_0^T \int_{\mathbb{R}_+} |\varphi(s,z)|^2 \nu(dz)\,ds < \infty \right).$$

Let us now give more specifications on the claim process

$$(2.2) \qquad S_t \stackrel{\triangle}{=} \int_0^{t+} \int_{\mathbb{R}_+} f(s,z) N_p(ds\,dz), \qquad t \geq 0.$$

We note that if the intensity $f(s,z) \equiv z$ and $\nu(\mathbb{R}_+) = \lambda > 0$, then $S_t$ is simply a compound Poisson process. Indeed, in this case one has $S_t = \sum_{0 \leq s < t, s \in D_p} \Delta S_t = \sum_{k \geq 1} \Delta S_{T_k} 1_{\{T_k \leq t\}}$, with $p_t \stackrel{\triangle}{=} \Delta S_t$ being a Poisson point process, $D_p \stackrel{\triangle}{=} \{t : p_t \neq 0\} = \bigcup_{k=1}^\infty \{T_k\}$, and $P\{p_{T_k} \in dz\} = \frac{1}{\lambda}\nu(dz)$, for all $k \geq 1$. Furthermore, $N_t \stackrel{\triangle}{=} \sum_{k=1}^\infty \mathbf{1}_{\{T_k \leq t\}}$ is a standard Poisson process with intensity $\lambda > 0$, and $S$ can be rewritten as $S_t = \sum_{k=1}^{N_t} p_{T_k}$ (cf. [8, 9]).

In what follows we shall make use of the following important assumption on the claim density $f$:

(H1) The random field $f \in F_p$, and it is continuous in $t$, and piecewise continuous in $z$. Furthermore, there exist constants $0 < d < L$ such that

$$(2.3) \qquad d \leq f(s,z,\omega) \leq L \qquad \forall (s,z) \in [0,\infty) \times \mathbb{R}_+, \ P\text{-a.s.}$$

We remark that the upper and lower bounds in (H1) reflects the simple fact in insurance: the *deductible* and *benefit limit*, and this is possible because $\nu(\mathbb{R}_+) < \infty$. Although mathematically we can replace such an assumption by certain integrability assumptions on both $f$ and $f^{-1}$ against the Lévy measure $\nu$, or that $f$ has a certain compact support in $z$, we prefer writing it in this simple way because of its practical meaning.

We now specify the *premium process* $\{c_t\}$ and the *safety loading* process $\{\rho_t\}$ in the reserve equation (1.1). In light of the well-known "*equivalence principle*" in actuarial mathematics (cf. Bowers et al. [2]), the premium process $\{c_t\}$ can be quantitatively characterized by the following equation:

$$(2.4) \qquad c_t = E\{\Delta S_t | \mathcal{F}_t^p\} = \int_{\mathbb{R}_+} f(t,z)\nu(dz) \qquad \forall t \geq 0, \ P\text{-a.s.}$$

Moreover, it is common to require that the premium and the expense loading satisfy the following "*net profit condition*":

$$(2.5) \qquad \operatorname*{essinf}_{\omega \in \Omega} \left\{ c_t(\omega)(1+\rho_t(\omega)) - \int_{\mathbb{R}_+} f(t,z,\omega)\nu(dz) \right\} > 0 \qquad \forall t \geq 0.$$

We therefore make the following standing assumption.

(H2) The safety loading process $\rho$ is a bounded, nonnegative $\mathbf{F}^p$-adapted process, and the premium process $c$ is an $\mathbf{F}^p$-adapted satisfying (2.4). Furthermore, the processes $c$, $\rho$, satisfy the "net profit condition" (2.5).

Note that if $f(t,z) \equiv z$, that is, the claim process is simply a compound Poisson, and $\rho$ is a constant, one has $c_s = c = \int_{\mathbb{R}_+} z\nu(dz) = \lambda E[U_1]$, where $U_1 = \Delta S_{T_1}$ is the jump size of the claim. In this case (2.5) becomes $c(1+\rho) > \lambda E(U_1)$, a usual net profit condition (cf. Asmussen–Nielsen [1]).

We now extend the reserve equation so that it contains the reinsurance and investment. We begin by the definition of a (generalized) reinsurance policy.

DEFINITION 2.1.  A (proportional) reinsurance policy is a random field $\alpha : [0, \infty) \times \mathbb{R}_+ \times \Omega \mapsto [0,1]$ such that $\alpha \in F_p$, and that for each fixed $z \in \mathbb{R}_+$, the process $\alpha(\cdot, z, \cdot)$ is predictable. Given a reinsurance policy $\alpha$, the part of the claim that a insurance company retains to itself during any time period $[t, t + \Delta t]$ is assumed to be $[\alpha * S]_t^{t+\Delta t}$, where $[\alpha * S]_0^t \stackrel{\triangle}{=} \int_0^t \int_{\mathbb{R}_+} \alpha(s, z) f(s, z) N_p(dz\, ds)$. In other words, the part of the claims it cedes to the reinsurer is $[(1-\alpha) * S]_t^{t+\Delta t}$.

We remark that the dependence of a reinsurance policy $\alpha$ on the spatial variable $z$ amounts to saying that the proportion can depend on the sizes of the claims, which is not unusual in practice. Although a traditional reinsurance policy as a predictable process $\alpha_t$, $t \geq 0$, might be simpler to treat from the modeling point of view, it is noted (as we shall see) that in general one may not be able to find an optimal strategy in such a form, unless $S_t$ has fixed size jumps [i.e., $\nu(dz)$ is a discrete measure]. But such a case is obviously not of significant interest because it will even exclude the compound Poisson claim processes.

The reserve equation with reinsurance can be argued heuristically using the so-called "*profit margin principle*" as follows. Let us denote the safety loading of the reinsurance company by $\rho^r$, and the modified safety loading of the original (cedent) company after reinsurance by $\rho^\alpha$. Consider an arbitrary small interval $[t, t+\Delta t]$, and denote $E_t^p\{\cdot\} = E\{\cdot|\mathcal{F}_t^p\}$. Then the following identity, which we call the "*profit margin principle*," should hold:

(2.6)
$$\underbrace{(1+\rho_t)E_t^p\{[1 * S]_t^{t+\Delta t}\}}_{\text{original premium}} - \underbrace{(1+\rho_t^r)E_t^p\{[(1-\alpha) * S]_t^{t+\Delta t}\}}_{\text{premium to the reinsurance company}}$$
$$= \underbrace{(1+\rho_t^\alpha)E_t^p\{[\alpha * S]_t^{t+\Delta t}\}}_{\text{modified premium}}.$$

Now assume that during this interval the reinsurance policy does not change in time. Using the assumption (H1) on $f$, one shows that, for any $\beta \in F_p$,

$$E_t\{[\beta * S]_t^{t+\Delta t}\} = \int_t^{t+\Delta t} \int_{\mathbb{R}_+} \beta(t,z) f(s,z) \nu(dz)\, ds$$

$$= \int_{\mathbb{R}_+} \beta(t,z) f(t,z) \nu(dz) \Delta t + o(\Delta t).$$

Now, approximating $E_t^p\{[\beta * S]_t^{t+\Delta t}\}$ by $\int_{\mathbb{R}_+} \beta(s,z) f(s,z) \nu(dz) \Delta t$ with $\beta = 1, \alpha, 1-\alpha$, respectively, in (2.6) and recalling (2.4), we obtain that

(2.7)
$$(1+\rho_t)c_t - (1+\rho_t^r) \int_{\mathbb{R}_+} (1-\alpha(t,z)) f(t,z) \nu(dz)$$
$$\approx (1+\rho_t^\alpha) \int_{\mathbb{R}_+} \alpha(t,z) f(t,z) \nu(dz).$$

Therefore, during $[t, t+\Delta t]$ the reserve changes as follows:

$$X_{t+\Delta t} - X_t$$

(2.8)
$$\approx c_t(1+\rho_t) \Delta t - (1+\rho_t^r) \int_{\mathbb{R}_+} (1-\alpha(t,z)) f(t,z) \nu(dz) \Delta t$$
$$- \int_t^{t+\Delta t} \int_{\mathbb{R}_+} \alpha(t,z) f(s,z) N_p(dz\,ds)$$
$$= (1+\rho_t^\alpha) \int_{\mathbb{R}_+} \alpha(t,z) f(s,z) \nu(dz) \Delta t$$
$$- \int_t^{t+\Delta t} \int_{\mathbb{R}_+} \alpha(t,z) f(s,z) N_p(dz\,ds).$$

For notational simplicity from now on let us denote

(2.9)
$$S_t^\alpha = \int_0^t \int_{\mathbb{R}_+} \alpha(t,z) f(s,z) N_p(dz\,ds),$$
$$m(t,\alpha) = \int_{\mathbb{R}_+} \alpha(t,z) f(s,z) \nu(dz).$$

Then (2.8) leads to the following equation for the reserve process:

(2.10)
$$X_t = x + \int_0^t (1+\rho_s^\alpha) m(s,\alpha)\,ds - S_t^\alpha$$
$$= x + \int_0^t (1+\rho_s^\alpha) m(s,\alpha)\,ds - \int_0^t \int_{\mathbb{R}_+} \alpha(s,z) f(s,z) N_p(ds\,dz).$$

REMARK 2.2. In the case when the reinsurance policy $\alpha$ is independent of $z$, we have $S_t^\alpha = \int_0^t \alpha(s) \int_{\mathbb{R}_+} f(s,z) N_p(dz\,ds) = \int_0^t \alpha(s)\,dS_s$ and $m(t,\alpha) = \alpha(s) \int_{\mathbb{R}_+} f(s,z) \nu(dz) = \alpha(s) c_s$, as we often see in the standard reinsurance framework. Also, we note that if $\rho^r = \rho^\alpha$ (hence equal to $\rho$!), then the reinsurance is called "*cheap*." But under the profit margin principle, we see that whether a reinsurance is cheap does not change the form of the reserve equation (2.10). From now on we shall drop the superscript "$\alpha$" from $\rho^a$ for simplicity, even when "noncheap" reinsurance is considered.

To conclude this section we now consider the scenario when an insurance company is allowed to invest part or all of its reserve in a financial market. We assume that the market has $k$ risky assets (stocks) and 1 riskless asset (bond or money market account). We model the dynamics of the market prices of the bond and stocks, denoted by $P_t^0$, $P_t^i$, respectively, where $t \geq 0$ and $i = 1, 2, \ldots, k$, which are described by the following stochastic differential equations:

$$
(2.11) \quad \begin{cases} dP_t^0 = r_t P_t^0 \, dt, \\ dP_t^i = P_t^i \left[ \mu_t^i \, dt + \sum_{j=1}^k \sigma_t^{ij} \, dW_t^j \right], \quad i = 1, 2, \ldots, k, \end{cases} \quad t \in [0, T],
$$

where $\{r_t\}$ is the interest rate, $\mu_t = (\mu_t^1, \mu_t^2, \ldots, \mu_t^k)^*$ is the appreciation rate, and $\sigma_t = (\sigma_t^{ij})_{i,j=1}^k$ is the volatility matrix. We shall make the following assumptions for the market parameters:

(H3) The processes $r$, $\mu$ and $\sigma$ are $\mathbf{F}^W$-adapted and bounded. Furthermore, the process $\sigma_t$ is uniformly nondegenerate, that is, $\sigma_t \sigma_t^* \geq \delta I$, $t \in [0, T]$, $P$-a.s., for some $\delta > 0$.

As usual, we assume that the market is liquid and the insurance company can trade continuously, and we denote the investment portfolio of the insurance company at each time $t$ by $\pi_t(\cdot) = (\pi_t^1, \ldots, \pi_t^k)$, where $\pi_t^i$ represents the fraction of its reserve $X_t$ allocated to the $i$th stock. Hence the amount of money that it puts into the $i$th stock would be $\pi_t^i X_t$, $i = 1, 2, \ldots, k$, and the rest of the money $X_t - \sum_{i=1}^k \pi_t^i X_t = (1 - \sum_{i=1}^k \pi_t^i) X_t$ will be put into the money market account. Furthermore, we denote the rate of the consumption of the insurance company to be an adapted process $D = \{D_t : t \geq 0\}$.

We should note that since we allow short selling and borrowing (with same interest rate). That is, we do not require that $\pi_t^i \geq 0$ and $\sum_{i=1}^k \pi_t^i \leq 1$, $\forall t \geq 0$. However, we do need some constraints on the portfolio process $\pi$ for technical reasons.

DEFINITION 2.3.  A portfolio process is an $\mathbb{R}^k$-valued, $\mathcal{F}$-adapted process $\pi$ such that

$$
(2.12) \qquad E\left\{ \int_0^T \|\pi_s X_s\|^2 \, ds \right\} = E\left\{ \sum_{i=1}^k \int_0^T |\pi_s^i X_s|^2 \, ds \right\} < \infty,
$$

where $X_t$ is the total reserve at time $t$. A consumption (rate) process is an $\mathcal{F}$-predictive nonnegative process $D$ satisfying

$$E\left\{\int_0^T D_s\,ds\right\} < \infty. \tag{2.13}$$

Following the idea of "self-financing," as before we now assume that during a small time duration $[t, t + \Delta t]$ the portfolio $\pi$, reinsurance policy $\alpha$ and the consumption rate $D$, as well as all the parameters are "freezed" at their values at time $t$, then it is easy to see that the reserve change during $[t, t + \Delta t]$ should be

$$\begin{aligned}
X_{t+\Delta t} \approx X_t &+ \underbrace{\sum_{i=1}^k \frac{\pi_t^i X_t}{P_t^i}\Delta P_t^i + \frac{(1-\sum_{i=1}^k \pi_t^i)X_t}{P_t^0}\Delta P_t^0}_{\text{investment gain}}\\
&+ \underbrace{(1+\rho_t)m(t,\alpha)\Delta t}_{\text{premium income}} - \underbrace{\int_t^{t+\Delta t}\int_{\mathbb{R}_+}\alpha(t,z)f(s,z)N_p(dz\,dt)}_{\text{claim}}\\
&- \underbrace{D_t\Delta t}_{\text{consumption}}.
\end{aligned} \tag{2.14}$$

Letting $\Delta t \to 0$, and using the price equations (2.11), we see that the reserve process $X$ should now follow the SDE

$$\begin{aligned}
X_t = x &+ \int_0^t \{X_s[r_s + \langle \pi_s, \mu_s - r_s \mathbf{1}\rangle] + (1+\rho_s)m(s,\alpha)\}\,ds\\
&+ \int_0^t X_s\langle \pi_s, \sigma_s\,dW_s\rangle\\
&- \int_0^{t+}\int_{\mathbb{R}_+}\alpha(s,z)f(s,z)N_p(ds\,dz) - \int_0^t D_s\,ds, \qquad t \in [0,T].
\end{aligned} \tag{2.15}$$

We often call a portfolio/reinsurance pair $(\pi,\alpha)$ is "$D$-financing" (see, e.g., [11]) if the risk reserve $X$ satisfy (2.15).

**3. Admissibility of strategies.** In this section we analyze some natural constraints on the investment and reinsurance strategies. We have already mentioned that the constraint $\alpha \in [0,1]$ is intrinsic in order to have a sensible reinsurance problem. Another special, fundamental constraint for an insurance company is that the reserve should (by government regulation) be aloft, that is, at any time $t \geq 0$, the reserve should satisfy $X_t^{x,\pi,\alpha,D} \geq C$ for some constant $C > 0$ at all time. Mathematically, one can always take $C = 0$ (or by changing $x$ to $x - C \geq 0$). We henceforth have the following definition of the "admissibility" condition.

DEFINITION 3.1. For any $x \geq 0$, a portfolio/reinsurance/consumption triplet $(\pi, \alpha, D)$ is called "admissible at $x$" if the risk reserve process satisfies

$$X_0^{x,\pi,\alpha,D} = x; \qquad X_t^{x,\pi,\alpha,D} \geq 0 \qquad \forall t \in [0,T], \ P\text{-a.s.}$$

We denote the totality of all strategies admissible at $x$ by $\mathcal{A}(x)$.

We observe that if $\alpha = 0$ and $D = 0$, then the reserve equation (2.15) becomes a homogeneous linear SDE, and $X_t^{x,\pi,0,0} > 0$ holds for all $t \geq 0$. But in general the admissibility of a given strategy becomes a more delicate issue, which we shall address in the sequel.

We begin by deriving a necessary condition for an admissible strategy. In light of the standard approach in finance, $\alpha \equiv 0$ in our case (see, e.g., Karatzas–Shreve [10, 11]), we denote the *risk premium* of the market by $\theta_t \stackrel{\triangle}{=} \sigma_t^{-1}(\mu_t - r_t \mathbf{1})$, and the *discount factor* by $\gamma_t \stackrel{\triangle}{=} \exp\{-\int_0^t r_s \, ds\}$, $t \geq 0$. Define

$$(3.1) \qquad W_t^0 \stackrel{\triangle}{=} W_t + \int_0^t \theta_s \, ds,$$

$$(3.2) \qquad Z_t \stackrel{\triangle}{=} \exp\left\{-\int_0^t \langle \theta_s, dW_s \rangle - \frac{1}{2} \int_0^t \|\theta_s\|^2 \, ds\right\},$$

$$(3.3) \qquad Y_t \stackrel{\triangle}{=} \exp\left\{\int_0^t \int_{\mathbb{R}_+} \ln(1+\rho_s) N_p(ds\, dz) - \nu(\mathbb{R}^+) \int_0^t \rho_s \, ds\right\}.$$

Finally, the so-called *state-price-density* process is defined as $H_t \stackrel{\triangle}{=} \gamma_t Y_t Z_t$.

We now give two lemmas concerning the Girsanov–Meyer transformations that will be useful in our discussion. Consider the following change of measures on the measurable space $(\Omega, \mathcal{F}_T)$:

$$(3.4) \qquad \begin{aligned} dQ_Z &= Z_T \, dP; \\ dQ &= Y_T \, dQ_Z = Y_T Z_T \, dP. \end{aligned}$$

Then by the Girsanov theorem (cf., e.g., [10]) we know that the process $W^0$ is a Brownian motion under measure $Q_Z$. The following lemma lists some less obvious consequences.

LEMMA 3.2. *Assume* (H2). *Then, under probability measure* $P$, *the process* $\{Y_t\}_{t \geq 0}$ *is a square-integrable martingale. Furthermore, define* $dQ_Y = Y_T \, dP$ *on* $\mathcal{F}_T$, *then:*

(i) *the process $Z$ is a square-integrable $Q_Y$-martingale;*

(ii) *for any reinsurance policy $\alpha$, the process*

$$N_t^\alpha \triangleq \int_0^t (1+\rho_s) m(s,\alpha) \, ds$$
(3.5)
$$- \int_0^{t+} \int_{\mathbb{R}_+} \alpha(s,z) f(s,z) N_p(ds\,dz)$$

*is a $Q_Y$-local martingale;*

(iii) *the process $ZN^\alpha$ is a $Q_Y$-local martingale.*

PROOF. Let $\xi_t = \int_0^t \int_{\mathbb{R}_+} \ln(1+\rho_s) N_p(ds\,dz) - \Lambda \int_0^t \rho_s \, ds$, where $\Lambda = \nu(\mathbb{R}_+)$. Then $Y_t = \exp\{\xi_t\}$. Applying Itô's formula (cf., e.g., [8]) we have

$$Y_t = 1 - \Lambda \int_0^t Y_s \rho_s \, ds$$

$$+ \int_0^{t+} \int_{\mathbb{R}_+} \{\exp\{\xi_{s-} + \ln(1+\rho_s)\} - \exp\{\xi_{s-}\}\} N_p(ds\,dz)$$

(3.6) $\qquad = 1 - \Lambda \int_0^t Y_s \rho_s \, ds$

$$+ \int_0^{t+} \int_{\mathbb{R}_+} Y_{s-}[\exp\{\ln(1+\rho_s)\} - 1] N_p(ds\,dz)$$

$$= 1 + \int_0^{t+} \int_{\mathbb{R}_+} Y_{s-} \rho_s \tilde{N}_p(ds\,dz),$$

where $\tilde{N}_p(ds\,dz) \triangleq N_p(ds\,dz) - \nu(dz) \, ds$ is the compensated Poisson random measure. That is, $Y$ is a local martingale. On the other hand, note that

$$Y_t^2 = \exp\left\{ \int_0^{t+} \int_{\mathbb{R}_+} 2\ln(1+\rho_s) N_p(ds\,dz) - \Lambda \int_0^t 2\rho_s \, ds \right\}$$

$$= \exp\left\{ \int_0^{t+} \int_{\mathbb{R}_+} \ln(1+\rho_s)^2 N_p(ds\,dz) - \Lambda \int_0^t [(1+\rho_s)^2 - 1] \, ds \right\} e^{\Lambda \int_0^t \rho_s^2 \, ds}$$

$$\triangleq \tilde{Y}_t e^{\Lambda \int_0^t \rho_s^2 \, ds},$$

where $\tilde{Y}$ is the same as $Y$ but with $\rho$ being replaced by $(1+\rho)^2 - 1$. Thus, repeating the previous arguments one shows that $\tilde{Y}$ satisfy the following stochastic differential equation

(3.7) $\qquad \tilde{Y}_t = 1 + \int_0^{t+} \int_{\mathbb{R}_+} \tilde{Y}_{s-}[(1+\rho_s)^2 - 1] \tilde{N}_p(ds\,dz),$

hence $\tilde{Y}$ is a local martingale as well. Since $\tilde{Y}$ is positive, it is a supermartingale. Therefore $E\tilde{Y}_t \leq E\tilde{Y}_0 = 1$ for all $t \geq 0$. The boundedness of $\rho$ then leads to that $EY_t^2 \leq E\tilde{Y}_t e^{\Lambda \int_0^t \rho_s \, ds} < \infty$. Thus $Y$ is indeed a true martingale.

Now consider processes $Z$ and $N^\alpha$ under the probability measure $Q_Y$. Since $\rho$ is $\mathbf{F}^p$-adapted, it is independent of $W$ (under $P$). Thus $Y$ and $Z$ are independent under $P$. Note that $Z$ satisfies the SDE:

$$(3.8) \qquad Z_t = 1 - \int_0^t \theta_s Z_s\, dW_s,$$

it is a square-integrable martingale under $P$, whence under $Q_Y$, proving (i).

To see (ii) we need only show that the process $N^\alpha Y$ is a $P$-local martingale for any reinsurance policy $\alpha$. Indeed, applying Itô's formula, noting that $Y$ satisfies the SDE (3.6), and recalling the definition of $m(\cdot,\alpha)$ (2.9), we have

$$(3.9) \qquad \begin{aligned} N_t^\alpha Y_t &= \int_0^{t+}\!\!\int_{\mathbb{R}^+} N_s^\alpha Y_{s-}\rho_s \tilde{N}_p(ds\,dz) \\ &\quad + \int_0^t\!\!\int_{\mathbb{R}^+} Y_{s-}(1+\rho_s) m(s,\alpha)\, ds \\ &\quad - \int_0^{t+}\!\!\int_{\mathbb{R}_+} Y_{s-}\alpha(s,z) f(s,z) N_p(ds\,dz) \\ &\quad - \int_0^t\!\!\int_{\mathbb{R}^+} \alpha(s,z) f(s,z) Y_{s-}\rho_s \nu(dz)\, ds \\ &= \int_0^{t+}\!\!\int_{\mathbb{R}^+} N_s^\alpha Y_{s-}\rho_s \tilde{N}_p(ds\,dz) \\ &\quad - \int_0^{t+}\!\!\int_{\mathbb{R}^+} Y_{s-}\alpha(s,z) f(s,z) \tilde{N}_p(ds\,dz). \end{aligned}$$

Thus $N^\alpha Y$ is $P$-local martingale.

Finally, (iii) follows from an easy application of Itô's formula. The proof is complete. □

A direct consequence of Lemma 3.2 is the following corollary.

COROLLARY 3.3. *Assume* (H2) *and* (H3). *The process $W^0$ is also a $Q$-Brownian motion, and $N^\alpha$ is a $Q$-local martingale. Consequently, $N^\alpha W^0$ is a $Q$-local martingale.*

PROOF. We first check $W^0$. Note that $W^0(t)$ is still a continuous process under $Q$, and for $0 \le s \le t$, one has

$$(3.10) \qquad \begin{aligned} E^Q\{W_t^0 - W_s^0 | \mathcal{F}_s\} &= \frac{1}{Y_s Z_s} E\{Y_T Z_T (W_t^0 - W_s^0) | \mathcal{F}_s\} \\ &= \frac{1}{Y_s} E\{Y_T | \mathcal{F}_s\} \frac{1}{Z_s} E\{Z_T (W_t^0 - W_s^0) | \mathcal{F}_s\} \\ &= E^{Q_Z}\{W_t^0 - W_s^0 | \mathcal{F}_s\} = 0. \end{aligned}$$

In the above the first equality is due to the Bayes rule (cf. [10]), the second equality is due to the independence of $Y$ and $Z$ and in the third equality we used the Bayes rule again, together with the facts that $Y$ is a $P$-martingale and $W^0$ is a $Q_Z$-Brownian motion. Similarly, one can show that

$$(3.11) \quad E^Q\{(W_t^0 - W_s^0)^2|\mathcal{F}_s\} = E^{Q_Z}\{(W_t^0 - W_s^0)^2|\mathcal{F}_s\} = I_d(t-s),$$

where $I_d$ is the $d \times d$ identity matrix. Applying Lévy's theorem we see that $W^0$ is a Brownian motion under $Q$.

To see that $N^\alpha$ is a $Q$-local martingale we must note that the reinsurance policy $\alpha$ is assumed to be **F**-adapted, hence $N^\alpha$ is neither independent of $Y$, nor of $Z$. We proceed with a slightly different argument. First notice that by an extra stopping if necessary, we may assume that $N^\alpha$ is bounded, whence a $Q_Y$-martingale by Lemma 3.2(ii). Also, in this case the conclusion (iii) of Lemma 3.2 can be strengthened to that $N^\alpha Z$ is a $Q_Y$-martingale as well. Bearing these in mind we apply Bayes rule again and use Lemma 3.2(i) to get

$$(3.12) \begin{aligned} & E^Q\{N_t^\alpha - N_s^\alpha|\mathcal{F}_s\} \\ &= \frac{1}{Z_s} E^{Q_Y}\{Z_T(N_t^\alpha - N_s^\alpha)|\mathcal{F}_s\} \\ &= \frac{1}{Z_s} E^{Q_Y}\{\{E^{Q_Y}\{Z_T|\mathcal{F}_t\}N_t^\alpha - E^{Q_Y}\{Z_T|\mathcal{F}_s\}N_s^\alpha\}|\mathcal{F}_s\} \\ &= \frac{1}{Z_s} E^{Q_Y}\{Z_t N_t^\alpha - Z_s N_s^\alpha|\mathcal{F}_s\} = 0. \end{aligned}$$

Thus $N^\alpha$ is a $Q$-martingale. The last claim is obvious. The proof is complete. □

The following necessary condition for the admissible triplet $(\pi, \alpha, D)$, also known as the "*budget constraint*," is now easy to derive.

THEOREM 3.4. *Assume* (H2) *and* (H3). *Then for any* $(\pi, \alpha, D) \in \mathcal{A}(x)$, *it holds that*

$$(3.13) \qquad E\left\{\int_0^T H_s D_s\, ds + H_T X_T^{x,\alpha,\pi,D}\right\} \leq x,$$

*where*

$$(3.14) \qquad H_t = \gamma_t Y_t Z_t, \qquad \gamma_t = \exp\left\{-\int_0^t r_s\, ds\right\}.$$

PROOF. For simplicity we denote $X = X^{x,\pi,\alpha,D}$. Recall the reserve equation (2.15) and rewrite it as

$$
\begin{aligned}
X_t &= x + \int_0^t \{r_s X_s + X_s \langle \pi_s, \sigma_s(\theta_s\,ds + dW_s)\rangle\} \\
&\quad + \int_0^{t+} (1+\rho_s) m(s,\alpha)\,ds \\
&\quad - \int_0^{t+} \int_{\mathbb{R}_+} \alpha(s,z) f(s,z) N_p(ds\,dz) - \int_0^t D_s\,ds \\
&= x + \int_0^t \{r_s X_s + X_s \langle \pi_s, \sigma_s\,dW_s^0\rangle\} + N_t^\alpha - \int_0^t D_s\,ds,
\end{aligned}
\tag{3.15}
$$

where $\theta_t = \sigma_t^{-1}(\mu_t - r_t)$ is the risk premium. Clearly, applying Itô's formula we can then write the discounted reserve, denoted by $\tilde{X}_t \stackrel{\triangle}{=} \gamma_t X_t$, $t \geq 0$, as follows:

$$
\tilde{X}_t = x + \int_0^t \tilde{X}_s \langle \pi_s, \sigma_s\,dW_s^0\rangle + \int_0^t \gamma_s\,dN_s^\alpha - \int_0^t \gamma_s D_s\,ds, \qquad t \geq 0. \tag{3.16}
$$

Therefore, under the probability measure $Q$ defined by (3.4), the process

$$
\tilde{X}_t + \int_0^t \gamma_s D_s\,ds = x + \int_0^t \tilde{X}_s \langle \pi, \sigma_s\,dW_s^0\rangle + \int_0^t \gamma_s\,dN_s^\alpha
$$

is a local martingale. Further, the admissibility of $(\pi, \alpha, D)$ implies that the left-hand side is a positive process, hence it is a supermartingale under $Q$. It follows that

$$
x \geq E^Q\left\{\tilde{X}_T + \int_0^T \gamma_s D_s\,ds\right\} = E\left\{H_T X_T + \int_0^T H_s D_s\,ds\right\}, \tag{3.17}
$$

proving the theorem. □

We remark that the budget constraint (3.13) takes the same form as those often seen in the pure finance models without claims (cf., e.g., [11]). The difference is that the "discounting" is accomplished by a different state-price-density process $H$.

**4. Wider-sense strategies and the auxiliary market.** In this and the next section we shall study the existence of admissible strategies, or more precisely, that the set of admissible strategies, $\mathcal{A}(x)$, is indeed nonempty for any initial endowment $x$. Since such existence results usually depend heavily on the *martingale representation theorem* (cf., e.g., [10]), or more generally the backward stochastic differential equation, the constraint on the reinsurance policy ($|\alpha| \leq 1$) causes fundamental difficulties in such argument. Our main idea is to first introduce a class of "*wider-sense strategies*," which essentially takes away the constraint, and then to find conditions under which a wider-sense strategy is indeed a member of $\mathcal{A}(x)$. We begin in this section with a detailed description of the wider-sense strategies as well as a wider-sense admissibility, and will discuss the nonemptiness of $\mathcal{A}(x)$ in Section 5.

DEFINITION 4.1. We say that a triplet of **F**-adapted processes $(\pi, \alpha, D)$ a wider-sense strategy if $\pi$ and $D$ satisfy (2.12) and (2.13), respectively, and $\alpha \in F_p^2$ [see (2.1)]. Moreover, we call the process $\alpha$ in a wider-sense strategy a *pseudo-reinsurance policy*.

The following lemma gives the existence of the wider-sense strategies.

LEMMA 4.2. *Assume* (H1)–(H3). *Then for any consumption process $D$ and any $\mathcal{F}_T$-measurable nonnegative random variable $B$ such that $E(B) > 0$ and*

$$E\left\{\int_0^T H_s D_s \, ds + H_T B\right\} = x, \tag{4.1}$$

*there exist a $D$-financing portfolio process $\pi$ and a pseudo-reinsurance policy $\alpha$, such that the solution $X^{x,\pi,\alpha,D}$ to the SDE (2.15) satisfies*

$$X_t^{x,\pi,\alpha,D} > 0 \qquad \forall 0 \le t \le T \quad and \quad X_T^{x,\pi,\alpha,D} = B, \qquad P\text{-}a.s.$$

PROOF. Let the consumption rate process $D$ be given. Consider the following *backward stochastic differential equation* (BSDE) on the probability space $(\Omega, \mathcal{F}, P)$:

$$X_t = B - \int_t^T \left\{r_s X_s + \langle \varphi_s, \theta_s \rangle - D_s + \rho_s \int_{\mathbb{R}_+} \psi(s,z)\nu(dz)\right\} ds \\ - \int_t^T \langle \varphi_s, dW_s \rangle + \int_t^T \int_{\mathbb{R}_+} \psi(s,z) \tilde{N}_p(ds\,dz). \tag{4.2}$$

Extending the results of BSDE with jumps by Situ [18] and using the martingale representation theorem involving random measures (cf. [9] or Lemma 2.3 in [19]), it can be shown that the BSDE (4.2) has a unique (**F**-adapted) solution $(X, \varphi, \psi)$ satisfying

$$E\int_0^T \left\{|X_s|^2 + |\varphi_s|^2 + \int_{\mathbb{R}_+} |\psi(s,z)|^2 \nu(dz)\right\} ds < \infty.$$

Let us define $\alpha(t,z) \stackrel{\triangle}{=} \frac{\psi(t,z)}{f(t,z)}$, for all $(t,z) \in [0,\infty) \times \mathbb{R}_+$, $P$-a.s. Then by (H1) we see that $\alpha \in F_p^2$, thus it is a pseudo-reinsurance policy. We claim that $X_t > 0$, for all $t \geq 0$, $P$-a.s. Indeed, note that

$$-\int_t^T \rho_s \int_{\mathbb{R}_+} \psi(s,z)\nu(dz)\,ds + \int_t^T \int_{\mathbb{R}_+} \psi(s,z)\tilde{N}_p(ds\,dz)$$

$$= -\int_t^T \rho_s \int_{\mathbb{R}_+} \alpha(s,z)f(s,z)\nu(dz)\,ds + \int_t^T \int_{\mathbb{R}_+} \alpha(s,z)f(s,z)\tilde{N}_p(ds\,dz)$$

$$= -\int_t^T (1+\rho_s)m(s,\alpha)\,ds + \int_t^T \int_{\mathbb{R}_+} \alpha(s,z)f(s,z)N_p(ds\,dz),$$

we see that (4.2) can be written as

$$X_t = B - \int_t^T \{r_s X_s + \langle \varphi_s, \theta_s \rangle - D_s + (1+\rho_s)m(s,\alpha)\}\,ds$$

(4.3)
$$- \int_t^T \langle \varphi_s, dW_s \rangle + \int_t^T \int_{\mathbb{R}_+} \alpha(s,z)f(s,z)N_p(ds\,dz)$$

$$= B - \int_t^T \{r_s X_s - D_s\}\,ds - \int_t^T \langle \varphi_s, dW_s^0 \rangle + N_T^\alpha - N_t^\alpha,$$

where $W^0$ and $N^\alpha$ are defined as before. Recall from Corollary 3.3 that $W^0$ is a $Q$-Brownian motion and $N^\alpha$ is a $Q$-local martingale. Let $\{\tau_n\}$ be a sequence of stopping times such that $\tau_n \uparrow \infty$ and for each $n$, $N_t^{\alpha,n} \triangleq N_{t\wedge\tau_n}^\alpha$, $t \geq 0$, is a martingale. Now for any $t \in [0,T]$, and any $n \geq 1$ we apply Itô's formula to get

$$\gamma_{T\wedge\tau_n} X_{T\wedge\tau_n} + \int_{t\wedge\tau_n}^{T\wedge\tau_n} \gamma_s D_s\,ds$$

$$= \gamma_{t\wedge\tau_n} X_{t\wedge\tau_n} + \int_{t\wedge\tau_n}^{T\wedge\tau_n} \gamma_s \langle \varphi_s, dW_s^0 \rangle - \int_{t\wedge\tau_n}^{T\wedge\tau_n} \gamma_s\,dN_s^{\alpha,n}.$$

Taking conditional expectations $E^Q\{\cdot|\mathcal{F}_{t\wedge\tau_n}\}$ on both sides and noting that the two stochastic integrals are all $Q$-martingales we obtain from the optional sampling theorem that

$$\gamma_{t\wedge\tau_n} X_{t\wedge\tau_n} = E^Q\bigg\{\gamma_{T\wedge\tau_n} X_{T\wedge\tau_n} + \int_{t\wedge\tau_n}^{T\wedge\tau_n} \gamma_s D_s\,ds\bigg|\mathcal{F}_{t\wedge\tau_n}\bigg\}.$$

Letting $n \to \infty$ and applying the monotone convergence theorem we then have

(4.4)
$$\gamma_t X_t = E^Q\bigg\{\gamma_T B + \int_t^T \gamma_s D_s\,ds\bigg|\mathcal{F}_t\bigg\} \geq E^Q\{\gamma_T B|\mathcal{F}_t\} > 0$$

$$\forall t \in [0,T],\ P\text{-a.s.},$$

since $E(B) > 0$ and $D$ is nonnegative by assumption. In other words, we have proved that $P\{X_t > 0, \forall t \geq 0; X_T = B\} = 1$. Let us now define

$$(4.5) \qquad \pi_t \triangleq (\sigma_t^T)^{-1} \frac{\varphi_t}{X_t}, \qquad t \in [0, T].$$

Then, we see that $\pi$ satisfies (2.12), thanks to (H3), and that (4.3) can now be written as

$$(4.6) \quad dX_t = \{r_t X_t - D_t\} dt + X_t \langle \pi_t, \sigma_t \, dW_t^0 \rangle - dN_t^\alpha, \qquad X_T = B.$$

Comparing (4.6) to the reserve equation (3.16) we see that $X = X^{x,\pi,\alpha,D}$ holds if we can show that $X_0 = x$. But setting $t = 0$ in (4.4) and using the assumption (4.1) we have

$$X_0 = E^Q \left\{ \gamma_T X_T + \int_0^T \gamma_s D_s \, ds \right\} = E \left\{ H_T X_T + \int_0^T H_s D_s \, ds \right\} = x.$$

This proves the lemma. $\square$

We remark that in Lemma 4.2 $\alpha$ is only a pseudo-reinsurance policy. In the rest of the section we will modify the wider-sense strategy obtained above to construct an admissible strategy. We shall follow the idea of the so-called "*duality method*" introduced by Cvitanic and Karatzas [4] to achieve this goal.

We begin by recalling the *support function* of $[0, 1]$ (see [4, 17]):

$$(4.7) \qquad \delta(x) \triangleq \delta(x|[0,1]) \triangleq \begin{cases} 0, & x \geq 0, \\ -x, & x < 0. \end{cases}$$

And we define a subspace of $F_p^2$:

$$(4.8) \quad \mathcal{D} \triangleq \left\{ v \in F_p^2 : \sup_{t \in [0,R]} \int_{\mathbb{R}^+} v(t,z) \nu(dz) < C_R, \ P\text{-a.s.}, \ \forall R > 0 \right\}.$$

Recall also the linear functional $m(\cdot, \cdot) : [0, T] \times F_p^2 \mapsto \mathbb{R}$ defined by

$$m(t, \alpha) = \int_{\mathbb{R}_+} \alpha(t, z) f(t, z) \nu(dz) \qquad \forall \alpha \in F_p^2.$$

Let $v \in \mathcal{D}$ be given. We consider a fictitious market in which the interest rate and appreciation rate are perturbed in such a way that the asset prices follow the SDE

$$(4.9) \quad \begin{cases} dP_t^{v,0} = P_t^{v,0} \{r_t + m(t, \delta(v))\} dt, \\ dP_t^{v,i} = P_t^{v,i} \left\{ (\mu_t^i + m(t, \delta(v))) dt + \sum_{j=1}^k \sigma_t^{ij} dW_t^j \right\}, \qquad i = 1, \ldots, k. \end{cases}$$

Next, recall the reserve equation (3.15) and definition (2.9) we have

$$X_t = x + \int_0^t r_t X_t \, dt + \int_0^t X_t \langle \pi_t, \sigma_t \, dW_t^0 \rangle$$

(4.10)
$$+ \int_0^t \int_{\mathbb{R}_+} (1 + \rho_s) \alpha(s,z) f(s,z) \nu(dz) \, ds$$

$$- \int_0^{t+} \int_{\mathbb{R}_+} \alpha(s,z) f(s,z) N_p(ds\, dz) - \int_0^t D_s \, ds.$$

Now corresponding to the auxiliary market we define a general (fictitious) expense loading function $\rho^v(s,z,x) \triangleq \rho_s + v(s,z)x$. Using such a loading function and repeating the previous argument one shows that the reserve equation (4.10) will become

$$X_t^v = x + \int_0^t X_s^v [r_s + m(s, \delta(v))] \, ds + \int_0^t X_s^v \langle \pi_s, \sigma_s \, dW_s^0 \rangle$$

$$+ \int_0^t \int_{\mathbb{R}^+} [1 + \rho_s + v(s,z) X_s^v] \alpha(s,z) f(s,z) \nu(dz) \, ds$$

$$- \int_0^{t+} \int_{\mathbb{R}^+} \alpha(s,z) f(s,z) N_p(ds\, dz) - \int_0^t D_s \, ds$$

(4.11)
$$= x + \int_0^t X_s^v [r_s + m(s, \alpha v + \delta(v))] \, ds + \int_0^t (1 + \rho_s) m(s, \alpha) \, ds$$

$$+ \int_0^t X_s^v \langle \pi_s, \sigma_s \, dW_s^0 \rangle - \int_0^{t+} \int_{\mathbb{R}^+} \alpha(s,z) f(s,z,\cdot) N_p(ds\, dz)$$

$$- \int_0^t D_s \, ds$$

$$= x + \int_0^t X_s^v r_s^{\alpha,v} \, ds + \int_0^t X_s^v \langle \pi_s, \sigma_s \, dW_s^0 \rangle + N_t^\alpha - \int_0^t D_s \, ds,$$

where $r_t^{\alpha,v} = r_t + m(t, \alpha v + \delta(v))$ could be thought of as a "fictitious" interest rate. We observe that for any pseudo-reinsurance strategy $\alpha$, the definition of $\delta(\cdot)$ implies that [suppressing variables $(t,z)$]:

(4.12) $\alpha v + \delta(v) = \alpha v \mathbf{1}_{\{v \geq 0\}} + (\alpha v - v) \mathbf{1}_{\{v < 0\}} = |v| \{\alpha \mathbf{1}_{\{v \geq 0\}} + (1 - \alpha) \mathbf{1}_{\{v < 0\}}\}$

and $r^{\alpha,v}$ reduces to the original interest rate if and only if $m(t, \alpha v + \delta(v)) = 0$. In particular, if $\alpha$ is a (true) reinsurance policy (hence $0 \leq \alpha \leq 1$), then it holds that

(4.13)
$$0 \leq \alpha(t,z) v(t,z) + \delta(v(t,z)) \leq |v(t,z)|$$

$$\forall (t,z) \in [0, \infty) \times \mathbb{R}_+, \ P\text{-a.s.}$$

The following modified wider-sense strategies will be useful in our future discussion.

DEFINITION 4.3. Let $v \in \mathcal{D}$. A wider-sense strategy $(\alpha, \pi, D)$ is called "$v$-admissible" if:

(i) $\int_0^T |m(t, av + \delta(v))| \, dt < \infty$, $P$-a.s.;
(ii) denoting $X^v = X^{v,x,\pi,\alpha,D}$, then $X_t^v \geq 0$, for all $0 \leq t \leq T$, $P$-a.s. We denote the totality of wider-sense $v$-admissible strategies by $\mathcal{A}^v(x)$.

We remark that if $v \in \mathcal{D}$ and $(\alpha, \pi, D) \in \mathcal{A}^v(x)$ such that

(4.14) $\begin{cases} 0 \leq \alpha(t,z) \leq 1; \\ \delta(v(t,z)) + \alpha(t,z)v(t,z) = 0, \end{cases} \quad dt \times \nu(dz)$-a.e., $P$-a.s.,

then $\alpha$ is a (true) reinsurance policy and $r_t^{\alpha,v} = r_t$, $t \geq 0$. Consequently, $X^v = X$ and $(\alpha, \pi, D) \in \mathcal{A}(x)$. To take a further look at the set $\mathcal{A}^v(x)$, let us define, for any $v \in \mathcal{D}$ and any $v$-admissible strategy $(\pi, \alpha, D)$,

(4.15) $\begin{cases} \gamma_t^{\alpha,v} \triangleq \exp\left\{-\int_0^t r_s^{\alpha,v} \, ds\right\} = \exp\left\{-\int_0^t [r_s + m(\alpha v + \delta(v))] \, ds\right\}, \\ H_t^{\alpha,v} \triangleq \gamma_t^{\alpha,v} Y_t Z_t, \quad t \in [0,T]. \end{cases}$

We have the following result.

PROPOSITION 4.4. *Assume* (H1)—(H3). *Then:*

(i) *for any $v \in \mathcal{D}$, and $(\pi, \alpha, D) \in \mathcal{A}^v(x)$, the following budget constraint still holds*

(4.16) $$E\left\{\int_0^T H_s^{\alpha,v} D_s \, ds + H_T^{\alpha,v} X_T^v\right\} \leq x;$$

(ii) *if $(\pi, \alpha, D) \in \mathcal{A}(x)$, then for any $v \in \mathcal{D}$ it holds that*

(4.17) $X^{v,x,\alpha,\pi,D}(t) \geq X^{x,\alpha,\pi,D}(t) \geq 0, \quad 0 \leq t \leq T, \; P\text{-a.s.}$

*In other words, $\mathcal{A}(x) \subseteq \mathcal{A}^v(x)$, $\forall v \in \mathcal{D}$.*

PROOF. (i) Recall from (3.3) and (3.2) the $P$-martingales $Y$ and $Z$, as well as the change of measure $dQ = Y_T Z_T \, dP$. Since $W^0$ is a $Q$-Brownian motion and $N^\alpha$ is a $Q$-local martingale, by the similar arguments as those in Theorem 3.4 one shows that

$$E\left\{\int_0^T H_s^{\alpha,v} D_s \, ds + H_T^{\alpha,v} X_T^v\right\} = E^Q\left\{\int_0^T \gamma_s^{\alpha,v} D_s \, ds + \gamma_T^{\alpha,v} X_T^v\right\} \leq x,$$

proving (i).

(ii) Let $x$ be fixed and let $(\alpha, \pi, D) \in \mathcal{A}(x)$. Since $\alpha$ is a (true) reinsurance policy, we have $\alpha(t,z) \in [0,1]$, for $dt \times \nu(dz)$-a.s. Thus $0 < \int_0^T m(s, \alpha v + \delta(v)) \, ds < \int_0^T |v(s,z)| \, ds < \infty$, thanks to (4.13). Thus denote $X = X^{x,\alpha,\pi,D}$, $X^v = X^{v,x,\alpha,\pi,D}$, and $\delta X \triangleq X^v - X$. We need only show that $\delta X_t \geq 0$, for all $t \geq 0$, a.s. Indeed, combining (2.15) and (4.11) we have

$$\delta X_t = \int_0^t (X_s^v r_s^{\alpha,v} - X_s r_s) \, ds + \int_0^t \delta X_s \langle \pi_s, \sigma_s \, dW_s^0 \rangle$$

(4.18)
$$= \int_0^t \delta X_s r_s^{\alpha,v} \, ds + \int_0^t \delta X_s \langle \pi_s, \sigma_s \, dW_s^0 \rangle$$

$$+ \int_0^t X_s m(s, \alpha v + \delta(v)) \, ds.$$

Viewing (4.18) as an linear SDE of $\delta X$ with $\delta X_0 = 0$, we derive from the *variation of parameter formula* that

$$\delta X_t = \mathcal{E}_t^{\alpha,v} \int_0^t [\mathcal{E}_s^{\alpha,v}]^{-1} X_s m(s, \alpha v + \delta(v)) \, ds,$$

where $\mathcal{E}^{\alpha,v} = \mathcal{E}(\xi^{\alpha,v})$ is the Doléans–Dade stochastic exponential of the semimartingale $\xi_t = \int_0^t r_s^{\alpha,v} \, ds + \int_0^t \langle \pi_s, \sigma_s \, dW_s^0 \rangle$, defined by

(4.19) $$\xi_t = \exp\left\{ \int_0^t \left[ r_s^{\alpha,v} + \frac{1}{2} |\pi_s \sigma_s|^2 \right] ds - \int_0^t \langle \pi_s, \sigma_s \, dW_s^0 \rangle \right\}$$

(cf., e.g., [10]). Note that $(\pi, \alpha, D) \in \mathcal{A}(x)$ implies that $X_t \geq 0$ for all $t \geq 0$ and that $m(t, \alpha v + \delta(v)) \geq 0$, thanks to (4.13). Consequently, $\delta X_t \geq 0$, for all $t \geq 0$, $P$-a.s. The proof is now complete. $\square$

**5. Existence of admissible strategies.** We are now ready to prove the existence of admissible strategies. Recall that Lemma 4.2 shows that the budget equation (4.1) implies the existence of a $D$-financing wider sense strategy. We will now look at the converse. To begin with, let us make some observations. Note that the perturbed reserve equation (4.11) can be rewritten as

$$X_t^v = x + \int_0^t \{[r_s + m(s, \delta(v) + \alpha v) + \langle \pi_s, \sigma_s \theta_s \rangle] X_s^v - D_s + \rho_s m(s, \alpha)\} \, ds$$

$$+ \int_0^t X_s^v \langle \pi_s, \sigma_s \, dW_s \rangle - \int_0^{t+} \int_{\mathbb{R}_+} \alpha(s,z) f(s,z) \tilde{N}_p(ds \, dz)$$

(5.1)
$$= x + \int_0^t \{[r_s + m(s, \delta(v))] X_s^v + m(s, \alpha v) X_s^v - D_s\} \, ds$$

$$+ \int_0^t X_s^v \langle \pi_s, \sigma_s \, dW_s^0 \rangle - \int_0^{t+} \int_{\mathbb{R}_+} \alpha(s,z) f(s,z) \tilde{N}_p^0(ds \, dz),$$

where $\tilde{N}_p^0(ds, dz) = \tilde{N}_p(ds, dz) - \rho_s \nu(dz)\, ds$. To simplify notation, we shall now denote, for any $v \in \mathcal{D}$ and $\eta \in F_p^2$,

$$(5.2) \qquad m^v(t, \eta) \triangleq \int_{\mathbb{R}_+} \eta(t, z) v(t, z) \nu(dz) \triangleq \overline{\eta}_t^v.$$

Thus, we have $m(t, \eta) = m^f(t, \eta)$, and we denote $m^1(t, \eta) = \overline{\eta}_t$. Let us now define

$$(5.3) \qquad \varphi_t^v = X_t^v \sigma_t^T \pi_t; \qquad \psi^v(t, z) = \alpha(t, z) f(t, z) \qquad \forall (t, z) \in [0, T] \times \mathbb{R}_+.$$

Then (5.1) becomes

$$(5.4) \qquad \begin{aligned} X_t^v = {} & x + \int_0^t \{[r_s + m(s, \delta(v))] X_s^v + \overline{\psi}_s^v X_s^v - D_s\} ds + \int_0^t \langle \varphi_s^v, dW_s^0 \rangle \\ & - \int_0^{t+} \int_{\mathbb{R}_+} \psi^v(s, z) \tilde{N}_p^0(ds\, dz). \end{aligned}$$

Recall that $W^0$ is a Brownian motion and $\tilde{N}^0$ is a compensated Poisson random measure under the probability measure $Q$, our analysis will depend heavily on the following BSDE deduced from (5.4): for any $B \in L^2(\Omega; \mathcal{F}_T)$, and $v \in F_p^2$,

$$(5.5) \qquad \begin{aligned} y_t = {} & B - \int_t^T \{r_s^v y_s + \overline{\psi}_s^v y_s - D_s\} ds - \int_t^T \langle \varphi_s, dW_s^0 \rangle \\ & + \int_t^T \int_{\mathbb{R}_+} \psi_s \tilde{N}_p^0(ds\, dz), \end{aligned}$$

where $r_t^v = r_t + m(t, \delta(v))$, $t \geq 0$. We should note here that this BSDE is not standard. Due to the multiplicative term $\overline{\psi}_s^v y_s$, which is neither Lipschitz nor linear growth in $(y^v, \psi^v)$. The following result is a special case of a general result by Liu [13], we provide a sketch of the proof that highlights the main idea.

LEMMA 5.1. *Assume* (H1)–(H3). *Assume further that processes $r$ and $D$ are all uniformly bounded. Then for any $v \in \mathcal{D}$ and $B \in L^\infty(\Omega; \mathcal{F}_T)$, such that $P\{B \geq 0\} = 1$, the BSDE (5.5) has a unique adapted solution $(y^v, \varphi^v, \psi^v)$.*

SKETCH OF THE PROOF. Denote $\kappa = e^{\|r^v\|_\infty} [\|B\|_\infty + T\|D\|_\infty]$, and define $\phi_\kappa : \mathbb{R} \to \mathbb{R}$ by $\phi_\kappa(y) = (y \wedge \kappa) \vee 0$. Consider the "truncated" version of (5.5):

$$(5.6) \qquad \begin{aligned} y_t^\kappa = {} & B - \int_t^T \{r_s^v y_s^\kappa + \overline{\psi}_s^{v,\kappa} \phi_\kappa(y_s^\kappa) - D_s\} ds - \int_t^T \langle \varphi_s^\kappa, dW_s^0 \rangle \\ & + \int_t^T \int_{\mathbb{R}_+} \psi_s^\kappa \tilde{N}_p^0(ds\, dz). \end{aligned}$$

Since this is now a BSDE with continuous, linear growth coefficients, following the techniques of [12] one shows that it has a solution $(y^\kappa, \varphi^\kappa, \psi^\kappa)$.

Next, define a new probability measure $Q^\kappa$ by $dQ^\kappa \triangleq \mathcal{E}^\kappa(T) \, dP$, where

$$(5.7) \quad \mathcal{E}^\kappa(t) \triangleq \exp\left\{ \int_0^t \int_{\mathbb{R}} \ln(1 + \phi_\kappa(y_s^\kappa)) N_p(ds\,dx) - \nu(\mathbb{R}) \int_0^t \phi_\kappa(y_s^\kappa) \, ds \right\}.$$

Then, under the measure $Q^\kappa$, BSDE (5.5) becomes

$$(5.8) \quad \begin{aligned} y_t^\kappa &= B - \int_t^T \{ r_s^v y_s^\kappa - D_s \} \, ds - \int_t^T \langle \varphi_s^\kappa, dW_s^\kappa \rangle \\ &\quad - \int_t^T \int_{\mathbb{R}} \psi_s^\kappa(x) \tilde{N}_p^\kappa(ds\,dx), \end{aligned}$$

where $\tilde{N}_p^\kappa(ds\,dx) \triangleq \tilde{N}_p^0(ds\,dx) - \phi_\kappa(y_s^\kappa)\nu(dx)\,ds$ is a $Q^\kappa$-martingale random measure. Then, a standard variation parameter argument shows that

$$(5.9) \quad y_t^\kappa = E^{Q^\kappa}\left\{ e^{-\int_t^T r_s^v \, ds} B + \int_t^T e^{-\int_t^u r_s^v \, ds} D_u \, du \Big| \mathcal{F}_t \right\}, \qquad Q^\kappa\text{-a.s.}$$

and thus $0 \leq |y_t^\kappa| \leq \kappa$, for all $t \in [0, T]$, $Q^\kappa$-a.s., and hence $\phi_\kappa(y_t^\kappa) = y_t^\kappa$, for all $t \in [0, T]$, $P$-a.s. To wit, $(y^\kappa, \varphi^\kappa, \psi^\kappa)$ satisfies the original BSDE (5.5). This proves the existence.

The uniqueness can also be proved using the boundedness of the solution and the change of measure technique; we omit it. □

We remark that for any $v \in F_p^2$, we can define a portfolio/pseudo-reinsurance policy pair from the solution $(y^v, \varphi^v, \psi^v)$ as

$$\pi_t^v = [\sigma_t^T]^{-1} \frac{\varphi_t^v}{y_t^v}; \qquad \alpha^v(t, z) = \frac{\psi^v(t, z)}{f(t, z)}.$$

Clearly, a necessary condition for $\alpha^v$ to be a true insurance policy is that

$$|\psi^v(t, z)| \leq |\alpha^v(t, z)||f(t, z)| \leq f(t, z) \leq L,$$

that is, $|\overline{\psi}_t^v| \leq L \sup_{t \in [0,T]} \int_{\mathbb{R}_+} |v(t, z)|\nu(dz) = L\|v\|_{\infty,\nu}$, where $L$ is the constant in (H1). In what follows we shall call the pair $(\pi^v, \alpha^v)$ the *portfolio/pseudo-reinsurance pair associated to $v$*.

With the help of Lemma 5.1, we now give a sufficient condition for the existence of the admissible strategy. The proof of this theorem borrows the idea of Theorem 9.1 in Cvitanic–Karatzas [4], modified to fit the current situation.

THEOREM 5.2. *Assume* (H1)–(H3). *Let $D$ be a bounded consumption process, and $B$ be any nonnegative, bounded $\mathcal{F}_T$-measurable random variable such that $E(B) > 0$. Suppose that for some $u^* \in \mathcal{D}$ whose associated portfolio/pseudo-reinsurance pair, denoted by $(\pi^*, \alpha^*)$, satisfies that*

$$E\left\{H_T^{\alpha^*,v}B + \int_0^T H_s^{\alpha^*,v}D_s\,ds\right\} \leq E\left\{H_T^{\alpha^*,u^*}B + \int_0^T H_s^{\alpha^*,u^*}D_s\,ds\right\} = x$$
$$\forall v \in \mathcal{D},$$

*where for any $v \in \mathcal{D}$,*

(5.10)
$$H_t^{\alpha^*,v} \triangleq \gamma_t^{\alpha^*,v}Y_tZ_t,$$
$$\gamma_t^{\alpha^*,v} \triangleq \exp\left\{-\int_0^t [r_s^v + m(s,\alpha^*v)]\,ds\right\}, \qquad t \geq 0.$$

*Then the triplet $(\pi^*, \alpha^*, D) \in \mathcal{A}(x)$. Further, the corresponding reserve $X^*$ satisfies $X_T^* = B$, $P$-a.s.*

PROOF. First note that with the given $u^* \in \mathcal{D}$ and the associated portfolio/pseudo-reinsurance pair $(\pi^*, \alpha^*)$, the BSDE (5.5) becomes

(5.11)
$$y_t^* = B - \int_t^T \{[r_s^{u^*} + m(s,\alpha^*u^*)]y_s^* - D_s\}\,ds - \int_t^T y_s^* \langle \pi_s^*, \sigma_s\,dW_s^0\rangle$$
$$+ \int_t^T \int_{\mathbb{R}_+} \alpha^*(s,z)f(s,z)\tilde{N}_p^0(ds\,dz).$$

Following the same arguments as that in Lemma 4.2 and using the assumption of the theorem, we can show that $y_T^* = B$, $y_t^* > 0$, for all $t \in [0,T]$, and $y_0^* = x$. In other words, we have shown that $(\pi^*, \alpha^*, D)$ is a $u^*$-admissible, wider-sense strategy. We shall prove that under the assumptions of the theorem, $\alpha^*$ is a true reinsurance policy and $y^* = X^{x,\pi^*,\alpha^*,D}$. But this is, according to the remark (4.14), amounts to showing that $0 \leq \alpha^*(t,z) \leq 1$, and $\delta(u^*(t,z)) + \alpha^*(t,z)u^*(t,z) = 0$, $dt \times \nu(dz) \times dP$-a.e.

In light of the arguments in [4], we introduce the following notation: for any $\lambda \in \mathcal{D}$, and $(t,z) \in [0,T] \times \mathbb{R}_+$, let

(5.12)
$$\delta^{u^*}(\lambda) \triangleq \begin{cases} -\delta(u^*), & \lambda = -u^*, \\ \delta(\lambda), & \text{otherwise.} \end{cases}$$

We then define, for $v \in \mathcal{D}$,

(5.13)
$$\begin{cases} x^*(v) \triangleq E\left[H_T^{\alpha^*,v}B + \int_0^T H_s^{\alpha^*,v}D_s\,ds\right]; \\ L_t^{*,v} \triangleq \int_0^t m(s,\delta^{u^*}(v-u^*))\,ds; \qquad M_t^{*,v} \triangleq \int_0^t m(s,\alpha^*(v-u^*))\,ds. \end{cases}$$

In particular, for any $\lambda \in \mathcal{D}$, we denote $L^\lambda \triangleq L^{*,u^*+\lambda}$ and $M^\lambda = M^{*,u^*+\lambda}$. We then define, for $\lambda \in \mathcal{D}$, a sequence of stopping times

$$(5.14) \quad \tau_n \triangleq \inf\left\{t \geq 0 : |L_t^\lambda| \vee |M_t^\lambda| \right.$$
$$\left. \vee \int_0^t \{[\gamma_s^{u^*} y_s^*(L_s^\lambda + M_s^\lambda)|\sigma_s^T \pi_s^*|]^2 + 1\}\,ds \geq n\right\} \wedge T.$$

It is clear from definition that $x^*(u^*) = y_0^* = x$ and $\tau_n \nearrow T$, as $n \to \infty$, $P$-a.s. for all $\lambda \in \mathcal{D}$.

Now for each $\lambda \in \mathcal{D}$, $\varepsilon \in (0,1)$ and $n \in \mathbb{N}$, we define a random field

$$(5.15) \quad u_{\varepsilon,n}^{*,\lambda}(t,z) \triangleq u^*(t,z) + \varepsilon \lambda(t,z) \mathbf{1}_{\{t \leq \tau_n\}}, \qquad (t,z) \in [0,T] \times \mathbb{R}_+.$$

Then clearly, $u_{\varepsilon,n}^{*,v} \in \mathcal{D}$. Furthermore, recalling (4.7) we see that $\delta$ is a convex function such that $\delta(cx) = c\delta(x)$ for all $c > 0$, and $x \in \mathbb{R}$. Thus it is easy to check that

$$(5.16) \quad \delta(u_{\varepsilon,n}^{*,\lambda}(t,z)) - \delta(u^*(t,z)) \leq \varepsilon \mathbf{1}_{\{t \leq \tau_n\}} \delta^{u^*}(\lambda(t,z))$$
$$\forall (t,z) \in [0,T] \times \mathbb{R}_+.$$

Now let $\lambda \in \mathcal{D}$ be fixed. For notational simplicity let us denote $H^* = H^{\alpha^*,u^*}$ and $H^{\varepsilon,n,\lambda} = H^{\alpha^*,u_{\varepsilon,n}^{*,\lambda}}$. Then, by the assumption of the theorem we have, for any $\varepsilon > 0$ and $n \in \mathbb{N}$,

$$0 \leq \frac{x^*(u^*) - x^*(u_{\varepsilon,n}^{*,\lambda})}{\varepsilon}$$
$$(5.17) \quad = \frac{1}{\varepsilon} E\left\{(H_T^* - H_T^{\varepsilon,n,\lambda})B + \int_0^T (H_s^* - H_s^{\varepsilon,n,\lambda})D_s\,ds\right\}$$
$$\triangleq E\{\Theta^{\varepsilon,n,\lambda}\},$$

where

$$(5.18) \quad \Theta^{\varepsilon,n,\lambda} \triangleq H_T^*\left(1 - \frac{H_T^{\varepsilon,n,\lambda}}{H_T^*}\right)\left(\frac{B}{\varepsilon}\right) + \int_0^T H_s^*\left(1 - \frac{H_s^{\varepsilon,n,\lambda}}{H_s^*}\right)\left(\frac{D_s}{\varepsilon}\right)ds.$$

We claim that, for each $n \in \mathbb{N}$,

$$(5.19) \quad \Theta^{\varepsilon,n,\lambda} \leq \kappa_n\left\{H_T^* B + \int_0^T H_s^* D_s\,ds\right\} \triangleq \Theta^n \qquad \forall \varepsilon > 0,$$

where $\kappa_n \triangleq \sup_{0 < \varepsilon < 1} \frac{1 - e^{-2\varepsilon n}}{\varepsilon}$. Indeed, recall from (5.10) and note (5.16) we see that

$$\frac{H_t^{\varepsilon,n,\lambda}}{H_t^*} = \frac{\gamma^{\alpha^*, u_{\varepsilon,n}^{*,\lambda}}}{\gamma^{\alpha^*,u^*}}$$

$$(5.20) \quad \begin{aligned} &= \exp\left\{-\int_0^t m(s, \delta(u_{\varepsilon,n}^{*;\lambda}) - \delta(u^*) + \alpha^*(u_{\varepsilon,n}^{*;\lambda} - u^*))\,ds\right\} \\ &\geq \exp\left\{-\varepsilon \int_0^{t\wedge\tau_n} m(s, \delta^{u^*}(\lambda) + \alpha^*\lambda)\,ds\right\} \\ &= \exp\{-\varepsilon(L_{t\wedge\tau_n}^\lambda + M_{t\wedge\tau_n}^\lambda)\} \geq e^{-2\varepsilon n}, \qquad t \in [0, T]. \end{aligned}$$

Thus (5.19) follows easily from (5.18). Furthermore, note that $E\Theta_n = \kappa y_0^* = \kappa_n x < \infty$, and $\varliminf_{\varepsilon\downarrow 0} \frac{1}{\varepsilon}\frac{H_t^{\varepsilon,n,\lambda}}{H_t^*} \geq -(L_{t\wedge\tau_n}^\lambda + M_{t\wedge\tau_n}^\lambda)$ for $t \in [0, T]$, we can apply Fatou's lemma to (5.17) to get

$$\begin{aligned} 0 &\leq \varlimsup_{\varepsilon\downarrow 0} \frac{x(u) - x(u_{\varepsilon,n}^v)}{\varepsilon} \leq E\left\{\varlimsup_{\varepsilon\downarrow 0} \Theta^{\varepsilon,n,\lambda}\right\} \\ &\leq E\left\{H_T^* B \varlimsup_{\varepsilon\downarrow 0} \frac{1}{\varepsilon}\left(1 - \frac{H_T^{\varepsilon,n,\lambda}}{H_T^*}\right) + \int_0^T H_s^* D_s \varlimsup_{\varepsilon\downarrow 0} \frac{1}{\varepsilon}\left(1 - \frac{H_s^{\varepsilon,n,\lambda}}{H_s^*}\right) ds\right\} \\ &\leq E\left\{H_T^* B(L_{\tau_n}^\lambda + M_{\tau_n}^\lambda) + \int_0^T H_s^* D_s(L_{s\wedge\tau_n}^\lambda + M_{s\wedge\tau_n}^\lambda)\,ds\right\} \\ &= E^Q\left\{\gamma_T^{u^*} B(L_{\tau_n}^\lambda + M_{\tau_n}^\lambda) + \int_0^T \gamma_s^{u^*} D_s(L_{s\wedge\tau_n}^\lambda + M_{s\wedge\tau_n}^\lambda)\,ds\right\}, \end{aligned}$$

where $Q$ is the probability measure defined as before. Now letting $n \to \infty$ and applying the dominated convergence theorem, we obtain that

$$(5.21) \quad 0 \leq E^Q\left\{\gamma_T^{u^*} B(L_T^\lambda + M_T^\lambda) + \int_0^T \gamma_s^{u^*} D_s(L_s^\lambda + M_s^\lambda)\,ds\right\}.$$

On the other hand, a simple application of Itô's formula to $\gamma^{u^*} y^*(L^\lambda + M^\lambda)$ from 0 to $\tau_n$ leads to that

$$(5.22) \quad \begin{aligned} \gamma_{\tau_n}^{u^*} y_{\tau_n}^*(L_{\tau_n}^\lambda + M_{\tau_n}^\lambda) &+ \int_0^{\tau_n}(L_s^\lambda + M_s^\lambda)\gamma_s^{u^*} D_s\,ds \\ &= \int_0^{\tau_n} \gamma_s^{u^*} y_s^* m(s, \delta(\lambda) + \alpha^*\lambda)\,ds + M_{\tau_n}, \end{aligned}$$

where

$$\begin{aligned} M_{t\wedge\tau_n} &\stackrel{\triangle}{=} \int_0^{t\wedge\tau_n}(L_s^\lambda + M_s^\lambda)\gamma_s^{u^*} y_s^*\langle \pi_s^*, \sigma_s\,dW_s^0\rangle \\ &\quad + \int_0^{t\wedge\tau_n} \gamma_s^{u^*}(L_s^\lambda + M_s^\lambda)\int_{\mathbb{R}_+} \alpha^*(s,z) f(s,z) N_p^0(ds\,dz) \end{aligned}$$

is a $Q$-martingale. First taking expectation on both sides of (5.22), then letting $n \to \infty$ and applying the dominated convergence theorem again we

obtain that

$$E^Q\left\{\gamma_T^{u^*} B(L_T^\lambda + M_T^\lambda) + \int_0^T (L_s^\lambda + M_s^\lambda)\gamma_s^{u^*} D_s \, ds\right\}$$

$$= E^Q\left\{\int_0^T \gamma_s^{u^*} y_s^* m(s, \delta(\lambda) + \alpha^*\lambda) \, ds\right\}.$$

Thus (5.21) becomes

(5.23) $$0 \leq E^Q\left\{\int_0^T \gamma_s^{u^*} y_s^* m(s, \delta^{u^*}(\lambda) + \alpha^*\lambda) \, ds\right\}.$$

Since $\lambda \in \mathcal{D}$ is arbitrary, we claim that this will lead to that

(5.24) $$\alpha^*(t,z,\omega)\lambda(t,z,\omega) + \delta^{u^*}(\lambda(t,z,\omega)) \geq 0,$$
$$dt \times d\nu \times dP\text{-a.e. } (t,z,\omega)$$

for all $\lambda \in \mathcal{D}$. Indeed, if not, we define $A \stackrel{\triangle}{=} \{(t,z,\omega) : \alpha^*(t,z,\omega)\lambda(t,z,\omega) + \delta^{u^*}(\lambda(t,z,\omega)) < 0\}$. Then it must hold that $Q\{\int_0^T \int_{\mathbb{R}_+} \mathbf{1}_A(t,z) \, dt \, \nu(dz) > 0\} > 0$. Since $\gamma_t^{u^*} y_t^* > 0$ for all $t \geq 0$, we have

$$I^\lambda(A) \stackrel{\triangle}{=} E^Q\left\{\int_0^T \gamma_s^{u^*} y_s^* m(s, \delta^{u^*}(\lambda \mathbf{1}_A) + \alpha^*\lambda\mathbf{1}_A) \, ds\right\} < 0.$$

Now for any $\eta > 0$ consider $\lambda^\eta \stackrel{\triangle}{=} \lambda\mathbf{1}_{A^c} + \eta\lambda\mathbf{1}_A \in \mathcal{D}$. Since $\eta$ is arbitrary, we can assume $\lambda^\eta \neq -u^*$ and thus $\delta^{u^*}(\lambda^\eta) = \delta(\lambda^\eta)$. Note again that $\delta(cx) = c\delta(x)$ for all $c > 0$ and $x \in \mathbb{R}$, it follows from (5.23) that

(5.25) $$0 \leq E^Q\left\{\int_0^T \gamma_s^{u^*} y_s^* m(s, \delta^{u^*}(\lambda^\eta) + \alpha^*\lambda^\eta) \, ds\right\} = I^\lambda(A^c) + \eta I^\lambda(A).$$

Since $I^\lambda(A) < 0$, we have a contradiction for $\eta > 0$ sufficiently large. This proves (5.24).

To complete the proof we note that the set $\mathcal{D}$ contains all constants. Thus (5.24) implies that for any $r \in \mathbb{R}$, $\alpha^*(t,z,\omega)r + \delta(r) \geq 0$, $dt \times d\nu \times dP$-a.e. Using the continuity of $\delta$ we can deduce easily that $\delta(r) + \alpha^*(t,z,\omega)r \geq 0$, $\forall r \in \mathbb{R}$, $dt \times d\nu \times dP$-a.e. Consequently we get from the definition of $\delta$ that

$$-\alpha^*(t,z,\omega)r \leq \delta(r) = \begin{cases} 0, & r \geq 0, \\ -r, & r < 0, \end{cases} \quad dt \times d\nu \times dP\text{-a.e.}$$

That is, $\alpha^*(t,z,\omega) \in [0,1]$, $dt \times d\nu \times dP$-a.e. Furthermore, note that (5.24) and definition of $\delta^{u*}(\cdot)$ imply that for both $\lambda = \pm u^*$ it holds that

$$\alpha^*(t,z)\lambda(t,z) + \delta^{u^*}(\lambda(t,z)) \geq 0, \qquad dt \times \nu(dz) \times dP\text{-a.s.}$$

That is, $\alpha^*(t,z)u^*(t,z) + \delta(u^*(t,z)) = 0$, $dt \times \nu(dz) \times dP$-a.e. The proof is complete. $\square$

**6. Utility optimization.** In this section we study a general utility optimization problem under our reserve model, by combining the results established in the previous sections and the results of [4, 5]. We begin by introducing some necessary notation.

DEFINITION 6.1. A function $U:[0,\infty) \mapsto [-\infty,\infty]$ is called a "utility function" if it enjoys the following properties:

(i) $U \in C^1((0,\infty))$, and $U'(x) > 0$, for all $x \in \mathbb{R}$;
(ii) $U'(\cdot)$ is strictly decreasing;
(iii) $U'(\infty) \triangleq \lim_{x\to\infty} U'(x) = 0$.

We denote $\text{dom}(U) \triangleq \{x \in [0,\infty); U(x) > -\infty\}$.

We should note that our definition of a utility function is slightly different from those in [4, 5] or [11] in that the $\text{dom}(U) \subseteq [0,\infty)$ instead of the whole real line. Recall (from [11], e.g.) that for any utility function, if we define $\bar{x} \triangleq \inf\{x \geq 0 : U(x) > -\infty\}$, and $U'(\bar{x}+) \triangleq \lim_{x \downarrow \bar{x}} U'(x)$. Then $U'(\bar{x}+) \in (0,+\infty]$. Furthermore, for each $x \in [0,\infty)$, let $I:(0,U'(\bar{x}+)) \mapsto (\bar{x},\infty)$ be the inverse of $U'$. Then $I$ is continuous and is strictly decreasing, and can be extended to $(0,\infty]$ by setting $I(y) = \bar{x}$ for $y \geq U'(\bar{x}+)$. Furthermore, one has

$$U'(I(y)) = \begin{cases} y, & 0 < y < U'(\bar{x}+), \\ U'(\bar{x}+), & U'(\bar{x}+) \leq y \leq \infty; \end{cases}$$

(6.1)
$$I(U'(x)) = x, \ \bar{x} < x < \infty.$$

In our optimization problem we consider the following "truncated version" of a utility function.

DEFINITION 6.2. A function $U$ is a "truncated utility function" if for some $K > 0$, $U$ is a utility function on $[0,K]$ but $U(x) = U(K)$ for all $x \geq K$. We call the interval $[0,K]$ the "effective domain" of $U$. Furthermore, we say that a truncated utility function is "good" if it satisfies $U'(\bar{x}+) < \infty$.

We note that for any utility function $U$, there exists a sequence of good truncated utility functions $\{U_n\}$ such that $U_n \to U$, as $n \to \infty$. Indeed, let $U$ be a utility function. For any $n > 0$, we define

(6.2)
$$\xi^n(x) = \begin{cases} \underline{x}^n - x + n, & 0 \leq x \leq \underline{x}^n, \\ U'(x), & \underline{x}^n \leq x \leq \bar{x}^n, \\ \dfrac{1}{n}, & x > \bar{x}^n, \end{cases}$$

where $\underline{x}^n$ and $\bar{x}^n$ are such that $U'(\underline{x}^n) = n$ and $U'(\bar{x}^n) = \frac{1}{n}$. Then it is fairly easy to check that the function defined by

$$U_n(x) = U(\underline{x}^n) - \frac{1}{2}(\underline{x}^n)^2 - n\underline{x}^n + \int_0^{x \wedge \bar{x}^n} \xi^n(y)\,dy, \qquad x \geq 0,$$

is a good truncated utility function with $\bar{x} = 0$, $U_n(0+) = U_n(0) = U(\underline{x}^n) - \frac{1}{2}(\underline{x}^n)^2 - n\underline{x}^n$, and $K = \bar{x}^n$. Clearly, $U_n(x) \to U(x)$, as $n \to \infty$, for all $x \in \mathbb{R}$.

We remark that if $U$ is a good truncated utility function, then we can assume without loss of generality that $\bar{x} = 0$, and $U'(0) < \infty$. Hence $U' : [0, K] \mapsto [U'(K-), U'(0)]$. If we define the inverse of $U'$ by

(6.3) $$I(y) = \inf\{x \geq 0 : U'(x) \leq y\}.$$

Then $I : [U'(K-), U'(0)] \mapsto [0, K]$ is continuous and strictly decreasing. We can also extend $I$ to $[0, \infty)$ by defining $I(y) = 0$ for $y \geq U'(0)$ and $I(y) = K$ for $y \in [0, U'(K-)]$. It is worth noting that such a function $I$ is then bounded over $[0, \infty)$ (!).

Now let $U$ be a good truncated utility function with effective domain $[0, K]$, and let us define the convex dual (or Legendre–Fenchel transform) of the function $U$ as follows:

(6.4) $$\tilde{U}(y) \stackrel{\triangle}{=} \max_{0 < x \leq K}\{U(x) - xy\}, \qquad 0 < y < \infty.$$

The following lemma shows that $\tilde{U}$ can be expressed in terms of $I$, just as the standard utility functions.

LEMMA 6.3. *Suppose that $U$ is a modified utility function, and let $\tilde{U}$ be its convex dual. Then it holds that*

(6.5) $$\tilde{U}(y) = U(I(y)) - yI(y) \qquad \forall y > 0.$$

PROOF. For each $y > 0$, consider the function $F(x) = U(x) - xy$. Differentiating with respect to $x$ we get $F'(x) = U'(x) - y$. If $y \in [U'(K-), U'(0)] = \text{Dom}(I)$, then we have $\tilde{U}(y) = \max_x F(x) = U(I(y)) - yI(y)$. If $y > U'(0)$, then we have $F'(x) = U'(x) - y < 0$ for all $x$, since $U'$ is decreasing. Thus $F(x)$ is decreasing and $\tilde{U}(y) = \max F(x) = F(0) = U(0)$. Since $I(y) = 0$ for all $y > U'(0)$, (6.3) still holds. Similarly, if $y \in [0, U'(K-))$, then $F'(x) = U'(x) - y > 0$ for all $x$. Thus $\tilde{U}(y) = \max F(x) = F(K) = F(I(y))$, for $y \in [0, U'(K-)]$. Thus (6.5) holds for all $y > 0$, proving the lemma. □

To formulate our optimization problem, we now consider a pair of functions $U_1 : [0, T] \times (0, \infty) \mapsto [-\infty, \infty)$ and $U_2 : [0, \infty) \mapsto [-\infty, \infty)$. As usual we denote the (partial) derivative of $U_1$ and $U_2$ with respect to $x$ by $U_1'$ and $U_2'$, respectively. We introduce the following definition that is based on the one in [11].

DEFINITION 6.4. A pair of functions $U_1:[0,T]\times(0,\infty)\mapsto[-\infty,\infty)$ and $U_2:[0,\infty)\mapsto[-\infty,\infty)$ is called a "(von Neumann–Morgenstern) preference structure" if:

(i) for each $t\in[0,T]$, $U^1(t,\cdot)$ is a utility function, such that the "subsistence consumption" defined by $\bar{x}_1(t)\triangleq\inf\{x\in\mathbb{R};U^1(t,x)>-\infty\}$ is continuous on $[0,T]$, and that both $U_1$ and $U_1'$ are continuous on the set $\mathcal{D}(U_1)\triangleq\{(t,x):x>\bar{x}^1(t),t\in[0,T]\}$;

(ii) $U_2$ is a utility function, with "subsistence terminal wealth" defined by $\bar{x}_2=\inf\{x:U_2'(x)>-\infty\}$.

Moreover, the pair $(U_1,U_2)$ is called a "modified preference structure" if $U_2$ is a good truncated utility function.

Our optimization problem is formulated as follows. Recall first the "risk neutral" measure $Q$ defined by (3.4). Let $(U_1,U_2)$, $x\in\mathbb{R}$ be a modified preference structure. We assume that the effective domain of $U_2$ is $[0,K]$. For any $(\pi,\alpha,D)\in\mathcal{A}(x)$, we define the total expected utility by

$$(6.6)\qquad J(x;\pi,\alpha,D)\triangleq E\bigg\{\int_0^T U_1(t,D_t)\,dt+U_2(X_T^{x,\alpha,\pi,D})\bigg\}.$$

The goal is to maximize $J(x;\pi,\alpha,D)$ over all $(\pi,\alpha,D)\in\mathcal{A}(x)$, and we denote the *value function* by

$$(6.7)\qquad V(x)\triangleq\sup_{(\pi,\alpha,D)\in\mathcal{A}(x)}J(x;\pi,\alpha,D).$$

We shall proceed along the lines of the duality method of [4, 5]. To be more precise, we shall first consider the fictitious market defined in Section 4 for find a candidate (wider-sense) optimal strategy, and then to verify that it is actually an optimal strategy in the strict sense, using Theorem 5.2.

To begin with, for given $v\in\mathcal{D}$ recall the set $\mathcal{A}_v(x)$, the $v$-admissible strategies, defined by Definition 4.3. For any $(\pi,\alpha,D)\in\mathcal{A}_v(x)$, the perturbed risk reserve, denoted by $X^v$ for simplicity, satisfies the following SDE under $Q$ [recall (4.11)]:

$$(6.8)\quad\begin{aligned}X^v(t)=x&+\int_0^t X_s^v[r_s+m(s,\delta(v)+\alpha v)]\,ds+\int_0^t X_s^v\langle\pi_s,\sigma_s\,dW_s^0\rangle\\&-\int_0^{t+}\int_{\mathbb{R}^+}\alpha(s,z)f(s,z)\tilde{N}_p^0(ds\,dz)-\int_0^t D_s\,ds,\qquad t\in[0,T].\end{aligned}$$

Let $\gamma^{\alpha,v}$ be the discounting factor defined by (4.15). Then as in Lemma 4.2 we can easily derive the following "fictitious" *budget constraint*:

$$(6.9)\quad x\geq E^Q\bigg\{\gamma_T^{\alpha,v}X_T^v+\int_0^T\gamma_s^{\alpha,v}D_s\,ds\bigg\}=E\bigg\{H_T^{\alpha,v}X_T^v+\int_0^T H_s^{\alpha,v}D_s\,ds\bigg\}.$$

Next, define

$$\mathcal{X}_v^\alpha(y) \stackrel{\triangle}{=} E\bigg\{ H_T^{\alpha,v} I_2(yH_T^{\alpha,v}) + \int_0^T H_t^{\alpha,v} I_1(t, yH_t^{\alpha,v})\, dt \bigg\},$$

(6.10)
$$0 < y < \infty,$$

where $I_1(t, \cdot)$ is the inverse of $U_1'(t, \cdot)$ and $I_2$ the inverse of $U_2'$. Since $I_1(t, \cdot)$ is continuous, strictly decreasing and $I_1(t, 0+) = \infty$, by the monotone convergence theorem we see that $\mathcal{X}_v^\alpha(0+) = \infty$ and $\mathcal{X}_v^\alpha(\cdot)$ is also continuous. On the other hand, note that $\lim_{x \to \infty} I_1(t, x) = \bar{x}^1(t)$ and $\lim_{x \to \infty} I_2(x) = \bar{x} = 0$, applying dominated convergence theorem we have

$$(6.11) \qquad \mathcal{X}_v^\alpha(\infty) \stackrel{\triangle}{=} \lim_{y \uparrow \infty} \mathcal{X}_v^\alpha(y) = E\bigg\{ \int_0^T H_t^{\alpha,v} \bar{x}^1(t)\, dt \bigg\}.$$

Furthermore, let $y_0 \stackrel{\triangle}{=} \sup\{y > 0; \mathcal{X}_v(y) > \mathcal{X}_v(\infty)\}$, then $y_0 \in (0, \infty]$ and one can show (as in [11]) that $X_v^\alpha(\cdot)$ is decreasing on the interval $(0, y_0)$. We can define the inverse of $\mathcal{X}_v^\alpha$ on $(0, y_0)$ by $\mathcal{Y}_v^\alpha(x) = \inf\{y : \mathcal{X}_v^\alpha(y) < x\}$. Then, $\mathcal{Y}_v^\alpha(x) \in (0, y_0)$, $\forall x \in (\mathcal{X}_v^\alpha(\infty), \infty)$.

Note that the definition of $\bar{x}_1(\cdot)$ tells us that $J(x; \pi, \alpha, D) > -\infty$ would imply

$$(6.12) \qquad D_t \geq \bar{x}_1(t), \qquad t \in [0, T],\ P\text{-a.s.}$$

But for all $(\alpha, \pi, D)$ satisfying (6.12), it holds that

$$(6.13) \quad E\bigg\{ H_T^{\alpha,v} X_T^v + \int_0^T H_t^{\alpha,v} D_t\, dt \bigg\} \geq E\bigg\{ \int_0^T H_t^{\alpha,v} D_t\, dt \bigg\} \geq \mathcal{X}_v^\alpha(\infty).$$

Therefore, the "fictitious" budget constraint (6.9) tells us that if $x < \mathcal{X}_v^\alpha(\infty)$, then (6.13) [whence (6.12)] cannot hold. Consequently one has $V(x) = -\infty$.

We now study the case when $\mathcal{X}_v^\alpha(\infty) < x$ holds. We note that this condition is rather unusual since it couples the initial state and the control. We now introduce a subset of $\mathcal{A}_v(x)$:

$$(6.14) \qquad \mathcal{A}_v'(x) \stackrel{\triangle}{=} \{(\pi, \alpha, D) \in \mathcal{A}_v(x) : x > \mathcal{X}_v^\alpha(\infty)\}.$$

Clearly, we need only consider the problem of maximizing

$$(6.15) \qquad \tilde{J}(D, B) \stackrel{\triangle}{=} E\bigg\{ \int_0^T U_1(t, D(t))\, dt + U_2(B) \bigg\}$$

over all pairs $(D, B)$, where $D$ is a consumption process and $B \in L^\infty_{\mathcal{F}_T}(\Omega)$, subject to the budget constraint

$$(6.16) \qquad E\bigg\{ \int_0^T H_t^{\alpha,v} D_t\, dt + H_T^{\alpha,v} B \bigg\} \leq x.$$

We again use the usual "Lagrange multiplier" method. For all $y > 0$, $v \in \mathcal{D}$ and $\alpha \in F_p^2$, let us try to maximize the following functional of $(D, B)$:

$$
\begin{aligned}
J_v^\alpha(D, B; x, y) &\triangleq E\left\{\int_0^T U_1(t, D(t))\, dt + U_2(B)\right\} \\
&\quad + y\left(x - E\left\{\int_0^T H_t^{\alpha,v} D_t\, dt + H_T^{\alpha,v} B\right\}\right) \\
&= xy + E\int_0^T [U_1(t, D(t)) - y H_t^{\alpha,v} D_t]\, dt \\
&\quad + E[U_2(B) - y H_T^{\alpha,v} B].
\end{aligned}
\tag{6.17}
$$

But recalling the definition of the convex duals of $U_1$ and $U_2$ we see that

$$
J_v^\alpha(D, B; x, y) \leq xy + E\left\{\int_0^T \tilde{U}_1(t, y H_t^{\alpha,v})\, dt + \tilde{U}_2(y H_T^{\alpha,v})\right\} \tag{6.18}
$$

with equality holds if and only if $D_t^{\alpha,v} = I_1(t, y H_t^{\alpha,v})$, $0 \leq t \leq T$, and $B^{\alpha,v} = I_2(y H_T^{\alpha,v})$, $P$-a.s.

We note that unlike the usual situations in finance (see, e.g., [4, 5]), the maximizer $D^{\alpha,v}$ and $B^{\alpha,v}$ depends on the reinsurance policy $\alpha$ as well. Since $B^{\alpha,v}$ is in a place of being the terminal reserve, our solution to the optimization problem therefore has a novel structure, which we now describe.

For any $y \in \mathbb{R}_+$ and $v \in \mathcal{D}$, consider the following so-called "forward–backward SDEs": for $t \in [0, T]$,

$$
\begin{cases}
H_t = 1 + \displaystyle\int_0^t H_s[r_s + m(s, \delta(v) + \alpha v)]\, ds - \int_0^t H_s \langle \theta_s, dW_s\rangle \\
\qquad + \displaystyle\int_0^t \int_{\mathbb{R}_+} H_{s-} \rho_s \tilde{N}_p(ds\, dz); \\
X_t = I_2(y H_T) - \displaystyle\int_t^T \{X_s[r_s + m(s, \delta(v) + \alpha v) + \langle \pi_s, \sigma_s \theta_s\rangle] \\
\hspace{6em} + (1 + \rho_s) m(s, \alpha)\}\, ds \\
\qquad - \displaystyle\int_t^T X_s \langle \pi_s, \sigma_s\, dW_s\rangle \\
\qquad + \displaystyle\int_t^T \int_{\mathbb{R}^+} \alpha(s, z) f(s, z) N_p(ds\, dz) + \int_t^T I_1(s, y H_s)\, ds.
\end{cases}
\tag{6.19}
$$

Denote the adapted solution to (6.19), if it exists, by $(H^{y,v}, X^{y,v}, \pi^{y,v}, \alpha^{y,v})$. Comparing the reserve equation (2.15) to the backward SDE in (6.19) we see that for a fixed $x$ and $D_t^{y,v} \triangleq I_1(t, y H_t^{y,v})$ for $t \geq 0$, the triplet $(\pi^{y,v}, \alpha^{y,v}, D^{y,v}) \in \mathcal{A}_v(x)$ if and only if $X_0^{y,v} = x$. But this is by no means clear from the FBSDE alone. In fact, we can only hope that for each $x$, there exists a $y = y(x)$, so that $X_0^{y(x),v} = x$. The following theorem is new.

THEOREM 6.5. *Assume* (H1)–(H3). *Let* $(U_1, U_2)$ *be a modified preference structure. The following two statements are equivalent:*

(i) *for any* $x \in \mathbb{R}$, $B^* \triangleq I_2(\mathcal{Y}(x)H_T)$ *and* $D_t^* \triangleq I_1(t, \mathcal{Y}(x)H_t)$, $t \geq 0$, *satisfy*

$$(6.20) \quad V(x) = E\left\{\int_0^T U_1(t, D_t^*)\,dt + U_2(B^*)\right\} = \sup_{(\pi,\alpha,D) \in \mathcal{A}(x)} J(x;\pi,\alpha,D),$$

*where* $\mathcal{Y}(x)$ *is such that*

$$(6.21) \quad x = E\left\{\int_0^T I_1(t, \mathcal{Y}(x)H_t)\,dt + I_2(\mathcal{Y}(x)H_T)\right\};$$

(ii) *there exists a* $u^* \in \mathcal{D}$, *such that the FBSDE (6.19) has an adapted solution* $(H^*, X^*, \pi^*, \alpha^*)$, *with $y$ satisfying*

$$(6.22) \quad x = E\left\{\int_0^T I_1(t, yH_t^*)\,dt + I_2(yH_T^*)\right\}.$$

*In particular, if* (i) *or* (ii) *holds, then* $(\pi^*, \alpha^*, D^*) \in \mathcal{A}(x)$ *is an optimal strategy for the utility maximization insurance/investment problem.*

PROOF. We first assume (i) and prove (ii). By assumption there exists a portfolio and reinsurance pair $(\pi^*, \alpha^*)$ such that $(\pi^*, \alpha^*, D^*) \in \mathcal{A}(x)$, $X_T^{\pi^*,\alpha^*,D^*} = B^*$, and that

$$J(x; \pi^*, \alpha^*, D^*) = V(x) = E\left\{\int_0^T U_1(t, D_t^*)\,dt + U_2(B^*)\right\}.$$

Since $\alpha^*(t, z) \in [0, 1]$, we can define a random field $u^*$ by

$$u^*(t, z) = \mathbf{1}_{\{\alpha^*(t,z)=0\}} - \mathbf{1}_{\{\alpha^*(t,z)=1\}} = \begin{cases} 1, & \alpha^*(t,z) = 0; \\ -1, & \alpha^*(t,z) = 1; \\ 0, & \text{otherwise,} \end{cases}$$

so that, by virtue of (4.12), $\delta(u^*) + \alpha^* u^* = |u^*|\{\alpha^* \mathbf{1}_{\{u^* \geq 0\}} + (1-\alpha^*)\mathbf{1}_{\{u^* < 0\}}\} \equiv 0$. Consequently, we must have $m(\cdot, \delta(u^*) + \alpha^* u^*) = 0$, $\gamma^{\alpha^*, u^*} = \gamma$, and $H^{\alpha^*, u^*} = H$. Note that the process $H$ satisfies the SDE

$$(6.23) \quad H_t = 1 + \int_0^t H_s r_s\,ds - \int_0^t H_s \langle \theta_s, dW_s \rangle + \int_0^t \int_{\mathbb{R}_+} H_{s-}\rho_s \tilde{N}_p(ds\,dz)$$

and the reserve $X^* \triangleq X^{\pi^*,\alpha^*,D^*}$ satisfies the SDE

$$(6.24) \quad \begin{aligned} X_t^* &= x + \int_0^t \{X_s^*[r_s + \langle \pi_s^*, \sigma_s\theta_s\rangle] + (1+\rho_s)m(s,\alpha^*)\}\,ds \\ &\quad + \int_0^t X_s^* \langle \pi_s^*, \sigma_s\,dW_s\rangle \\ &\quad - \int_0^t \int_{\mathbb{R}_+} \alpha^*(s,z)f(s,z)N_p(ds\,dz) - \int_0^t D_s^*\,ds \end{aligned}$$

and the terminal condition $X_T^* = B^* = I_2(\mathcal{Y}(x)H_T)$. Combining (6.23) and (6.24) we see that $(H, X^*, \pi^*, \alpha^*)$ actually satisfies the FBSDE (6.19) with $y = \mathcal{Y}(x)$ and $v = u^*$. Note that when $H^* = H$ the equations (6.21) and (6.22) coincide, we proved the statement (ii).

We now assume (ii) holds and try to prove (i). Suppose that for some $u^* \in \mathcal{D}$, FBSDE (6.19) has an adapted solution $(H^*, X^*, \pi^*, \alpha^*)$ with $y = \mathcal{Y}_{u^*}^{\alpha^*}(x) \triangleq \mathcal{Y}^*(x)$ [i.e., $y$ satisfies (6.22)]. We define

$$D_t^* = I_1(t, \mathcal{Y}^*(x)H_t^*), \qquad t \geq 0, \qquad B^* \triangleq I_2(\mathcal{Y}^*(x)H_T^*).$$

Since we have already seen that $(D^*, B^*)$ is a maximizer of the Lagrange multiplier problem for $J_{u^*}^{\alpha^*}(x, \mathcal{Y}^*(x); D, B)$, defined by (6.17), we must have

$$(6.25) \qquad x = E\left\{ H^* B^* + \int_0^T H_t^* D_t^* \, dt \right\}$$

and

$$(6.26) \qquad \begin{aligned} V^*(x) &= \sup_{(D,B)} J_{u^*}^{\alpha^*}(D, B; x, \mathcal{Y}_{u^*}^{\alpha^*}(x)) \\ &= E\left\{ \int_0^T U_1(t, D_t^*) \, dt + U_2(B^*) \right\}. \end{aligned}$$

Since $I_2$ is bounded by $K > 0$, and $H^*$ satisfy an SDE, we have $P\{B^* > 0\} = P\{\mathcal{Y}^*(x)H_T^* \leq U_2'(0)\} > 0$. Namely, $B^*$ is bounded and $E(B^*) > 0$. Furthermore, for any other $v \in \mathcal{D}$, the budget constraint (6.16) tells us that

$$E\left\{ \int_0^T H_t^{\alpha^*,v} D_t^* \, dt + H_T^{\alpha^*,v} B^* \right\} \leq x = E\left\{ \int_0^T H_t^* D_t^* \, dt + H_T^* B^* \right\}.$$

Thus, applying Theorem 5.2 we can conclude that $(\alpha^*, \pi^*, D^*) \in \mathcal{A}(x)$, to wit,

$$\begin{cases} 0 \leq \alpha^*(t, z) \leq 1; \\ \delta(u^*(t, z)) + \alpha^*(t, z)u^*(t, z) = 0, \end{cases} \qquad (t, z) \in [0, T] \times \mathbb{R}_+, \ P\text{-a.s.},$$

and $X_0^* = x$. Thus, we must again have $H^* = H$, $\mathcal{Y}^*(x) = \mathcal{Y}(x)$, and $X^* = X^{x,\pi^*,\alpha^*,D^*}$. Consequently, we see that $(D^*, B^*)$ become the same as that defined in (i), and $V^*(x) = V(x) = E\{\int_0^T U_1(t, D_t^*) \, dt + U_2(B^*)\}$, proving (i).

The last claim is clear from the proof. This completes the proof. $\square$

We note that the statement (ii) in Theorem 6.5 is indeed more than the solvability of the FBSDE (6.19), due to the special choice of $y$. In fact in a sense $y$ itself becomes a part of the solution. This, along with the solvability of FBSDE (6.19), forms a class of new problems in the theory of BSDEs, and it is currently under investigation. We hope to be able to address this issue in our future publications.

Wachovia Corporation
CP18 North Carolina 1001
201 South College Street
Charlotte, North Carolina 28244
USA
E-mail: yuping.liu@wachovia.com

Department of Mathematics
University of Southern California
Los Angeles, California 90089
USA
E-mail: jinma@usc.edu